\documentclass{amsart}

\usepackage{amssymb}

\newtheorem{theorem}{Theorem}[section]
\newtheorem{lemma}[theorem]{Lemma}
\newtheorem{corollary}[theorem]{Corollary}
\newtheorem{proposition}[theorem]{Proposition}

\newtheorem{remark}[theorem]{Remark}

\newenvironment{Proof}{\removelastskip\par\medskip
\noindent{\em Proof.} \rm}{\penalty-20\null\hfill$\square$\par\medbreak}

\numberwithin{equation}{section}
\newcommand{\Ai}{\text{Ai\,}}

\newcommand{\tr}{\text{tr\,}}

\begin{document}
\setcounter{page}{1}

\title[Non-intersecting paths, random tilings and random matrices]
{Non-intersecting paths, random tilings and random matrices}
\author[K.~Johansson]{Kurt Johansson}

\address{
Department of Mathematics,
Royal Institute of Technology,
S-100 44 Stockholm, Sweden}

\email{kurtj@math.kth.se}

\begin{abstract}
We investigate certain measures induced by families of
non-\linebreak intersecting paths in domino tilings of the Aztec diamond, rhombus
tilings of an abc-hexagon, a dimer model on a cylindrical brick
lattice and a growth model. The measures obtained, e.g. the Krawtchouk
and Hahn ensembles, have the same structure as the eigenvalue measures
in random matrix theory like GUE, which can in fact can be obtained
from non-intersecting Brownian motions. The derivations of the
measures are based on the Karlin-McGregor or
Lindstr\"om-Gessel-Viennot method. We use the measures to show
some asymptotic results for the models.
\end{abstract}

\maketitle

\section{Introduction}
We begin by summarizing some facts about the Gaussian Unitary Ensemble
of random hermitian matrices. We will see analogues of these in the
random tiling problems discussed below.
The {\it Gaussian Unitary Ensemble} (GUE) is the probability measure
\begin{equation}\label{1.1}
\frac 1{\mathcal{Z}_N}e^{-\tr M^2}dM
\end{equation}
on the space of all $N\times N$ hermitian matrices, which is
isomorphic to $\mathbb{R}^{N^2}$, and $dM$ is the Lebesgue measure on
this space. The induced measure on the $N$ real eigenvalues
$\lambda_1,\dots, \lambda_N$ of $M$ is given by, \cite{Me},
\begin{equation}\label{1.2}
\phi_{N,\text{GUE\,}}(\lambda)d^N\lambda=
\frac 1{Z_N}\Delta_N(\lambda)^2\prod_{j=1}^Ne^{-\lambda_j^2}d^N\lambda
\end{equation}
where
\begin{equation}\label{1.3}
\Delta_N(\lambda)=\det(\lambda_j^{N-k})_{j,k=1}^N=
\prod_{1\le j<k\le N}(x_j-x_k).
\end{equation}
The probability of finding $m$ eigenvalues in infinitesimal intervals
$d\lambda_1,\dots,d\lambda_m$ around $\lambda_1,\dots,\lambda_m$ is
given by
$R_{m,N}(\lambda_1,\dots,\lambda_m)d\lambda_1,\dots,d\lambda_m$, where
$R_{m,N}$ is the {\it m-point correlation function} 
\begin{equation}\label{1.4}
R_{m,N}(\lambda_1,\dots,\lambda_m)=\frac
{N!}{(N-m)!}\int_{\mathbb{R}^{N-m}}\phi_{N,\text{GUE\,}}(
\lambda_1,\dots,\lambda_N)d\lambda_{m+1},\dots,d\lambda_N.
\end{equation}
The correlation functions are given by determinants,
\begin{equation}\label{1.5}
R_{m,N}(\lambda_1,\dots,\lambda_m)=\det(K_N(\lambda_i,\lambda_j))_{i,j=1}^m,
\end{equation}
where
\begin{equation}\label{1.6}
K_N(x,y)=\frac{\kappa_{N-1}}{\kappa_N}\frac{h_N(x)h_{N-1}(y)-h_{N-1}(x)h_N(y)}
{x-y} \left(e^{-x^2-y^2}\right)^{1/2},
\end{equation}
and $h_n(x)=\kappa_n x^n+\dots$ are the orthonormal Hermite
polynomials. 

We can think of the eigenvalues as points on $\mathbb{R}$. Using the
formulas above and asymptotics for Hermite polynomials we can obtain a
limiting random point process with correlation functions
\begin{equation}\label{1.7}
R_m(x_1,\dots,x_m)=\det(\frac{\sin\pi(x_i-x_j)}{\pi(x_i-x_j)})_{i,j=1}^m.
\end{equation}
Hence we obtain a determinantal random point process given by the sine
kernel, \cite{So2}. To get this limit we rescale so that the mean
distance between the eigenvalues (points) equals 1 and then use,
\begin{equation}\label{1.8}
\lim_{N\to\infty}\frac 1{\sqrt{2N}\rho(u)}K_N\left(\sqrt{\frac N2}u+\frac
{\xi}{\sqrt{2N}\rho(u)},\sqrt{\frac N2}u+\frac
{\eta}{\sqrt{2N}\rho(u)}\right)=\frac{\sin\pi(\xi-\eta)}{\pi(\xi-\eta)},
\end{equation}
where $\rho(u)=\frac 1{2\pi}\sqrt{(4-u^2)_+}$ is the {\it Wigner
semicircle law} which describes the asymptotic density of the
eigenvalues. 

Let us recall two limit theorems for the fluctuations of the
eigenvalues. Denote by $\#(u,v)$ the number of eigenvalues in the
interval $[u\sqrt{N/2},v\sqrt{N/2}]$, $u<v$, $|u|,|v|<2$. Then,
 
\begin{equation}\label{1.9}
\frac{\#(u,v)-\mathbb{E}_N[\#(u,v)]}{\sqrt{\frac 2{\pi^2}\log N}}
\end{equation}
converges in distribution to the standard normal, \cite{CL},
\cite{So1}, \cite{So2}. If $\lambda_{\text{max}}$\linebreak
 $=\max_{1\le j\le
  N}\lambda_j$ is the largest eigenvalue, then
\begin{equation}\label{1.10}
\mathbb{P}_N\left[\frac{\lambda_{\text{max}}
-\sqrt{2N}}{\sqrt{2}N^{-1/6}}\le\xi\right]\to F(\xi) 
\end{equation}
as $N\to\infty$, $\xi\in\mathbb{R}$, where
\begin{equation}\label{1.11}
F(\xi)=\det(I-A)_{L^2(\xi,\infty)}.
\end{equation}
Here $A$ is the operator on $L^2(\xi,\infty)$ with kernel (the {\it
  Airy kernel})
\begin{equation}\label{1.12}
A(x,y)=\frac{\Ai(x)\Ai'(y)-\Ai'(x)\Ai(y)}{x-y},
\end{equation}
\cite{TW1}. The distribution function (\ref{1.11}) is called the {\it
  Tracy-Widom distribution}. It follows from (1.5) and the Fredholm
expansion that 
\begin{align}\label{1.13}
\mathbb{E}_N\left[\prod_{j=1}^N(1+g(\lambda_j))\right]&=\sum_{n=0}^N\frac
1{n!} \int_{\mathbb{R}^n}\det(K_N(x_i,x_j)g(x_j))_{i,j=1}^nd^nx\\
&=\det(I+K_Ng)_{L^2(\mathbb{R})}.\notag
\end{align}
If we take $g(x)=-\chi_{(t,\infty)}(x)$ in (\ref{1.13}) and use the
asymptotics of the Hermite polynomials close to the largest zero, we
can prove (\ref{1.10}).

The GUE eigenvalue measure, (\ref{1.2}), was obtained above from the
measure (\ref{1.1}) on random hermitian matrices. We will now show
how we can obtain (\ref{1.2}) in a completely different way using
non-intersecting Brownian motions. This type of problem has been
studied under the name of vicious walkers or domain walls in the
statistical physics litterature, see 
\cite{Fi}, \cite{Fo1}, \cite{Fo2}, \cite{Fo2'}
and references in these papers. Consider $N$ 1-dimensional Brownian
motions starting at the points $0,\dots, N-1$ at time 0 and ending at
the same points $0,\dots, N-1$ at time $2T$. Let
$p_{N,T}(x_1,\dots,x_N)$ denote the probability density that at time
$T$ the particles are at the positions $x_1<\dots<x_N$ conditioned on
the event that the paths have not intersected in the whole time
interval $[0,2T]$. If $p_t(x,y)=(2\pi t)^{-1/2}\exp (-(x-y)^2/2t)$ is
the transition kernel for Brownian motion, then, by a theorem of Karlin
and McGregor, \cite{KM}, see also \cite{Jo4},
\begin{align}\label{1.14}
p_{N,T}(x_1,\dots,x_N)&=\frac 1{Z_N}\det(p_T(j-1,x_k))_{j,k=1}^N
\det(p_T(x_j,k-1))_{j,k=1}^N\\
&=\frac 1{Z_N}\left(\det(p_T(j-1,x_k))_{j,k=1}^N\right)^2,\notag
\end{align}
where
\begin{equation}\label{1.15}
Z_N=\frac
1{N!}\int_{\mathbb{R}^N}\left(\det(p_T(j-1,x_k))_{j,k=1}^N\right)^2d^Nx.
\end{equation}
Note that, because of symmetry, we can consider (\ref{1.14}) as a
probability measure on $\mathbb{R}^N$ and remove the $N!$ in
(\ref{1.15}). It follows from \cite{TW2} that the probability measure
on $\mathbb{R}^N$ with density (\ref{1.14}) has determinantal
correlation functions analogous to (\ref{1.5}) but with a different
kernel, see \cite{Jo4}. We can now obtain GUE as follows. Let
$T\to\infty$ and rescale the $x_j$:s by $\sqrt{T}$ so that they do not
move away to infinity. Then,
\begin{equation}\label{1.16}
\lim_{T\to\infty}
p_{N,T}(x_1\sqrt{T},\dots,x_N\sqrt{T})=\phi_{N,\text{GUE\,}}(x). 
\end{equation}
To show this note that,
\begin{equation}
\det(p_T(j-1,x_k))_{j,k=1}^N=\frac 1{(2\pi
  t)^{-N/2}}\prod_{j=1}^Ne^{((j-1)^2+x_j^2)/2T}\det(e^{(j-1)x_k/2T})_{j,k=1}^N
\end{equation}
and use the formula for a Vandermonde determinant. The choice of
initial and final positions made above is not necessary for the result
but simplifies the computations. Thus, we need not look upon
(\ref{1.2}) as something which necessarily comes from random matrices.
For another relation between random matrices and Brownian motion see
\cite{Bar}, \cite{GTW}.

The present paper can be seen as a continuation of 
the paper \cite{Jo2} where several examples of analogues of GUE on a
discrete space was given. In these ensembles we have analogues of the
results for GUE discussed above but we obtain the so called discrete
sine kernel, in the limit instead of the ordinary sine
kernel. The common theme of the present paper is non-intersecting
paths, which are discrete analogues of the non-intersecting Brownian
motions just mentioned, and thus it is reasonable to expect that we
will have features analogous to those of GUE. In section 2 we will
consider random tilings of the Aztec diamond of size $n$, 
\cite{EKLP}, which can
be described using certain non-intersecting paths.
By using so called zig-zag paths in the tiling, we obtain the
Krawtchouk ensemble, \cite{Jo2}, which can be used to 
analyze several properties
of the random tiling. By a result in \cite{JPS} 
the shape of the so called temperate region in
a random tiling is closely related to the corner growth model in
\cite{Jo1}, which is a generalization of the longest increasing
subsequence problem for random permutations. 
This gives a new approach to the asymptotic fluctuation 
results in \cite{BDJ} and
\cite{Jo1} involving the Tracy-Widom distribution. 
We will also study the fluctuations of the domino height
function which describes the diling. It has Gaussian fluctuations with
variance of order $\log N$, a fact that is related to the Gaussian
fluctuations in the number of eigenvalues in an interval in a GUE
matrix. This type of result has been conjectured in \cite{PrSt}.
We will also discuss the relation between the equilibrium measure for
the Krawtchouk ensemble and the arctic ellipse.
The corner growth model can be
generalized further and this leads to the so called Schur measure
introduced on \cite{Ok1}. The Schur measure can also be analyzed using
non-intersecting paths as will be demonstrated in section 3.
We will also, in section 4, 
consider rhombus tilings of a hexagon, \cite{CLP}, which
are related to boxed plane partitions. 
These tilings can also be described by certain non-intersecting random
walk paths and the intersection of these paths with a fixed line leads
to the so called Hahn ensemble.
In this problem we will not
compute the detailed asymptotics, but we will discuss the equilibrium
measure for the Hahn ensemble and its relation to the arctic ellipse 
phenomenon. Finally, in section 5, we will
analyze certain aspects of a dimer model on a brick (hexagonal) lattice on a
cylinder. 
Here we also have non-intersecting paths but the number of paths is
not fixed like in the other examples.
The methods used are
very close to the arguments used to compute the correlation functions
in \cite{TW2}

There are many papers in the statistical physics litterature related
to the present paper, e.g. \cite{BEO}, \cite{Fi}, \cite{Fo1}, \cite{Fo2},
\cite{Fo2'}, \cite{GOV}, \cite{K}, \cite{KZ-J}, 
\cite{YNS} and \cite{Z-J} 
Connections between random permutations and
the so called random turns model, which gives certain non-intersecting
paths has beeen discussed in \cite{Bai}, \cite{Fo3},
and \cite{Fo4}. Other relevant papers are \cite{Fisch}, \cite{FK1}, 
\cite{FK2}, \cite{Gr}, \cite{HW} and
\cite{Pi}.
\section{The Aztec diamond}
\subsection{Basic definitions}
The {\it Aztec diamond}, $A_n$, of order $n$
is the union of all lattice squares $[m,m+1]\times[l,l+1]$,
$m,l\in\mathbb{Z}$,
that lie inside the region $\{(x,y)\,;\, |x|+|y|\le n+1\}$. A {\it
  domino}  is a closed $1\times 2$ or $2\times 1$ rectangle in
$\mathbb{R}^2$ with corners in $\mathbb{Z}^2$, and a {\it tiling} of a
region $R\subseteq \mathbb{R}^2$ by dominoes is a set of dominoes
whose interiors are disjoint and whose union is $R$. Let
$\mathcal{T}(A_n)$ denote the set of all domino tilings of the Aztec
diamond.

We can equivalently think of a tiling as a dimer
configuration. Consider the graph $G$ with vertices at $(\frac
12,\frac 12)+\mathbb{Z}^2$ and edges between nearest neighbour
vertices. A {\it dimer} is simply an edge in $G$, and if the edge goes
between the verices $v_1$ and $v_2$ we say that the dimer {\it covers} 
$v_1$ and $v_2$. Let $G_n$ be the subgraph of $G$ where all vertices
lie in $A_n$. A {\it dimer configuration} in $G_n$ is a set of dimers
in $G_n$ such that all vertices are covered by exactly one dimer. This
is clearly equivalent to a tiling of $A_n$ via the identification: a
dimer between $v_1$ and $v_2$ corresponds to a domino covering the two
lattice squares with centers $v_1$ and $v_2$.

 Colour the Aztec diamond in a checkerboard fashion so that
the leftmost square in each row in the top half is white. A horizontal domino
is {\it north-going} (N) if its leftmost square is white, otherwise it is 
{\it south-going} (S). Similarly, a vertical domino is {\it
  west-going} (W) if its 
upper square is white, otherwise it is {\it east-going} (E). Two dominoes are
\it adjacent \rm 
if they share an edge, and a domino is adjacent to the boundary
if it shares an edge with the bundary of the Aztec diamond. The 
\it north polar region \rm is defined to be the union of those north-going 
dominoes that are connected to the boundary by a sequence of adjacent 
north-going dominoes. The south, west and east polar regions are defined 
analogously. In this way a domino tiling partitions the Aztec diamond into 
four polar regions, where we have a regular brick wall pattern, and a fifth
central region, the \it temperate zone\rm, where the tiling pattern is 
irregular.

We will now define a one-to-one mapping from $\mathcal{T}(A_n)$ to
families of $n$ non-intersecting lattice paths. Consider an S-domino
which we place with corners at $(0,0), (2,0), (2,1), (0,1)$. Draw a
straight line from $(0,1/2)$ to $(2, 1/2)$. In this way we get a piece
of a path, and we do this for all S-dominoes. Similarly we can put a
W-domino so that it has corners at $(0,0), (1,0), (1,2), (0,2)$, and
then draw a straight line segment from $(0,1/2)$ to
$(1,3/2)$. Finally, on an E-domino placed at the same position we draw
a straight line from $(0,3/2)$ to $(1,1.2)$. Do this for all W- and
E-dominoes. We do not draw any line on an N-domino. Given a domino
tiling of $A_n$ we draw lines on the dominoes as just described. We
claim that this gives $n$ non-intersecting paths starting at
$A_j=(n+1-j,\frac 12-j)$ and ending at $E_j=(-n-1+j,\frac 12-j)$,
$1\le j\le n$. We call these paths {\it DR-paths of type I} (after
D. Randall, \cite{St}, p.277). To prove the claim we argue as
follows. Consider the black lattice square to the right of $E_j$. It
has to be covered by a W- or am S-domino. In both cases a path will
end at $E_j$ From the checkerboard colouring we see that we nust
obtain connected paths. Similarly, if we consider the white lattice
square to the left of $A_j$ it can only be covered by an E- or an
S-domino. Hence a path must start at the point $A_j$. Clearly, by
construction, the paths are non-intersecting and the claim is proved.

A convenient coordinate system for describing the paths is what we
call {\it coordinate system I} (CS-I). As origin we take $(n+1,1/2)$ and as
basis vectors $\mathbf{e}_I=(-1,-1)$, $\mathbf{f}_I=(-1,1)$. 
Let $\mathcal{L}_I$ be the
integer lattice in CS-I. The type I DR-paths are walks in
$\mathcal{L}_I$. They take steps $(1,0)$, $(0,1)$ or $(1,1)$ and they
have starting points $(k,0)$ and endpoints $(n+1,k)$, $1\le k\le
n$. Thus we obtain a map from domino tilings of $A_n$  to families of
$n$ nonintersecting type I DR-paths in $\mathcal{L}_I$ with the
specified initial and final positions. This map is a bijection. To see
this, fill in with dominoes along the paths using the marked
tiles. This is possible since the paths do not intersect. If we have a
white lattice square that is covered by a domino, then the black
square to the right is also empty, since otherwise the paths would not
be connected. Similarly, if a black square is not covered by a domino,
the white square to the left is not covered either. Hence, the squares
that are not already covered can be covered by N-dominoes. Clearly
this gives an inverse.

We can also define {\it type II DR-paths} which are complementary to the
type I paths. In this case the S-dominoes are unmarked, whereas the
N-dominoes have a horizontal segment in the middle. Furthermore we
interchange the marking on the W- and E-dominoes. In this way we
obtain paths from $A_j=(-n-1+j, j-1/2)$ to $E_j=(n+1-j, j-1/2)$, $1\le
j\le n$. In coordinate system II (CS-II) which has origin
$(-n-1,-1/2)$ and basis vectors $\mathbf{e}_{II}=(1,1),
\mathbf{f}_{II}
=(1,-1)$, we
obtain the same type of lattice paths as before.

Call the top type I DR-path the {\it level-1 path}. It is clear from the
definitions above that the north polar zone is exactly the part of the
Aztec diamond above the level-1 path, i.e. all dominoes, which have to
be N-domonioes, that lie above this path.

Let $\tau\in\mathcal{T}(A_n)$ be a tiling of the Aztec diamond and let
$v(\tau)$ denote the number of vertical dominoes in $\tau$. We define
a probability measure on $\mathcal{T}(A_n)$ by letting the horizontal
dominoes have weight 1 and the vertical dominoes weight $w$. Thus,
\begin{equation}\label{2.1}
\mathbb{P}[\tau]=\frac{w^{v(\tau)}}{\sum_{\tau\in\mathcal{T}(A_n)}
w^{v(\tau)}}.
\end{equation}
If we take $w=1$ we obtain the uniform distribution on
$\mathcal{T}(A_n)$. This can alternatively be viewed as a probability
measure on the DR-paths, where we put the weight 1 on the steps
$(1,1)$ (in CS-I or CS-II), which correspond to horizontal dominoes,
and the weight $w$ on the steps (1,0) or (0,1), which correspond to
vertical dominoes. The weight of $n$ given non-intersecting DR-paths
is the product of the weights on all steps and equals $w^{v(\tau)}$,
if $\tau$ is the tiling defined by the paths. The weight of a set of
non-intersecting DR-paths is the sum of the weights of all the
elements in the set.

Next, we will define another type of paths, the so called zig-zag
 paths, \cite{EKLP}, in the Aztec diamond. Consider the sequence of
white squares with opposite corners $Q_k^r=(-r+k, n+1-k-r)$,
$k=0,\dots, n+1$, where $r$, $1\le r\le n$, is fixed. A { \it zig-zag
 path} $Z_r$ in $A_n$ is a path of edges going around these white
 squares. When going from $Q_k^r$ to $Q_{k+1}^r$ we can go either
 first one step east and then one step south (an ES-step), or first
 one step south and then one step east (an SE-step). A domino tiling
 $\tau\in\mathcal{T}(A_n)$ defines a unique zig-zag path $Z_r(\tau)$
 from $Q_0^r$ to $Q_{n+1}^r$ if we require that the zig-zag path does
 not intersect the dominoes. There will be exactly $r$ ES-steps, and
 hence $n+1-r$ SE-steps along the zig-zag path. This can be proved
 using the domino height function defined below.

We associate the point $(r, n-k)$ in CS-I with the step
$Q_k^rQ_{k+1}^r$. Suppose that we have ES-steps at the points
$(r,h_j)$, $1\le j\le r$, in CS-I, i.e. $Q_{n-h_j}^rQ_{n-h_j+1}^r$ are
ES-steps. Then the zig-zag path is mapped one-to-one to $(h_1,\dots,
h_r)$, where $0\le h_1<\dots<h_r\le n$. This specifies the {\it zig-zag
(particle) configuration} $(h_1,\dots,h_r)$; we write
$p(Z_r)=(h_1,\dots,h_r)$. We can also associate the step
$Q_k^rQ_{k+1}^r$ with the point $(n+1-r,k)$ in CS-II, and we will then
have SE-steps at the points $(n+1-r, n-k_j)$, $1\le j\le n+1-r$, where
$k_1<\dots<k_{n+1-r}$. We call $(k_1,\dots, k_{n+1-r})$ the {\it zig
  zag (hole) configuration}, and write $h(Z_r)=(k_1,\dots,
k_{n+1-r})$. The next lemma gives the relation between the DR-paths
and the zig-zag paths.

\begin{lemma}\label{lem2.1}
The points $(r, h_j)$ in CS-I are the {\it last}
  positions on $x_I=r$ of the type I DR-paths starting at $(k,0)$,
  $1\le k\le r$ in CS-I. Similarly, $(n+1-r, n-k_j)$ are the {\it
  last} positions on $x_{II}=n+1-r$ in CS-II of the type II DR-paths
  starting at $(k,0)$, $1\le k\le n+1-r$. Also,
\begin{equation}\label{2.2}
\{h_1,\dots,h_r\}\cup\{k_1,\dots,k_{n+1-r}\}=\{0,\dots,n\}.
\end{equation}
\end{lemma}
 
\begin{Proof}
If we have an ES-step around a white square, then this square is
covered by an S- or a W-domino. In both cases it follows, from the
definition of the type I DR-paths, which are walks in the integer
lattice $\mathcal{L}_I$ in CS-I, that the DR-path must intersect the
S-step in the ES-step, and after that go to a point with a larger
$x_I$-coordinate. Similarly, if we have an SE-step around a white
square, then this square is covered by an N- or an E-domino, in which
case a type I DR-path does not intersect neither the S- nor the
E-step. The proof of the second statement in the lemma is analogous,
and (\ref{2.2}) follows from the definition of the particle and hole
configurations, and the definition of CS-I and CS-II.
\end{Proof}

If $p[Z_r(\tau)]=(h_1^r,\dots,h_r^r)$, $1\le r\le n$, then the
position of the rightmost particle, $h_r^r$, describes the level-1
type I DR-path. As noted above, the region above this DR-path is the
north polar zone, and hence we can investigate the shape of the north
polar zone using the positions of the rightmost particles. We will
return to this in sects. 2.3 and 2.4 below.

Let us recall the definition of the {\it (domino) height function}
associated with a given tiling, \cite{EKLP}. Let $u$ and $v$ be two
adjacent lattice points (vertices), in the basic coordinate system,
such that the edge connecting them is not covered by a domino. If the
edge from $u$ to $v$ has a black square to its left, $h(v)=h(u)+1$,
and if it has a white square to its left, $h(v)=h(u)-1$. 
Note that the value of the height function is uniquely
determined apart from an overall additive constant. 
We can fix it by requiring $h(n,0)=0$.
If $u$ and $v$ are
two adjacent lattice points, then $|h(u)-h(v)|=3$ if the edge is
covered by a domino, otherwise $|h(u)-h(v)|=1$. From this it follows
that $h(Q_0^r)=2n-(2r-1)$, $h(Q_{n+1}^r)=2r-1$ and
\begin{equation}\label{2.3}
h(Q_k^r)-h(Q_{k+1}^r)=
\begin{cases}
-2, &\text{if ES-step}\\
2,  &\text{if SE-step.}
\end{cases}
\end{equation}
Consequently, we can use the zig-zag configurations to determine the
height at a given point. From (\ref{2.3}) it follows that there are
exactly $r$ ES-steps in $Z_r(\tau)$.

\subsection{The Krawtckouk ensemble}
The {\it Krawtchouk ensemble}, \cite{Jo2}, is a probability measure on
$\{0,\dots, K\}^N$ defined by
\begin{equation}\label{2.4}
\mathbb{P}_{\text{Kr\,},N,K,p}[h]=\frac
1{Z_{N,K,p}}\Delta_N^2(h)\prod_{j=1}^N \binom{K}{h_j}p^{h_j}q^{K-h_j},
\end{equation}
where $0<p<1$, $q=1-p$, $1\le N\le K$, $h=(h_1,\dots,
h_N)\in\{0,\dots, K\}^N$ and
\begin{equation}\label{2.5}
Z_{N,K,p}=N!\left(\prod_{j=0}^{N-1}
\frac{j!}{(K-j)!}\right)K!^N(pq)^{N(N-1)/2}. 
\end{equation}
Set $w(x)=\binom{K}{x}p^xq^{K-x}$, $0\le x\le K$ and let
$\{p_k(x)\}_{k=0}^K$ be the normalized orthogonal polynomials with
respect to $w(x)$ on $\{0,\dots, K\}$, i.e.
\begin{equation}\label{2.5'}
\sum_{x=0}^K p_j(x)p_k(x)w(x)=\delta_{jk}.
\end{equation}
These are multiples of the ordinary Krawtchouk polynomials,
\cite{NSU}, and have the integral representation
\begin{equation}\label{2.6}
p_n(x)=\binom{K}{x}^{-1/2}(pq)^{-n/2}\frac 1{2\pi i}\int_\gamma\frac
{(1+qz)^x(1-pz)^{K-x}}{z^n} \frac{dz}z,
\end{equation}
where $\gamma$ is a circle centered at the origin with radius
$\le\min(1/p,1/q)$. The measure (\ref{2.4}) has determinantal
correlation functions, \cite{Me}, \cite{TW2},
\begin{equation}\label{2.6'}
\det(K_{\text{Kr\,},N,K,p}(x_i,x_j))_{i,j=1}^m
\end{equation}
where the {\it Krawtchouk kernel} is given by
\begin{align}\label{2.7}
K_{\text{Kr\,},N,K,p}(x,y)&=\sum_{n=0}^{N-1}p_n(x)p_n(y)(w(x)w(y))^{1/2}\\
&=\frac{\kappa_{N-1}}{\kappa_N}\frac{p_N(x)p_{N-1}(y)-p_{N-1}(x)p_N(y)}{x-y}
(w(x)w(y))^{1/2} \notag
\end{align}
where $\kappa_n=(n!)^{-1}\binom{K}{n}^{-1/2}(pq)^{-n/2}$ is the
leading coefficient in $p_n(x)$.

We will now prove, using the DR-paths, that the measure on the zig-zag
configurations induced by the probability measure (\ref{2.1}) on the
tilings is the Krawtchouk ensemble. In the special case of uniform
distribution on the set of tilings this was proved in \cite{Jo2} using
formulas from \cite{EKLP}.
\begin{theorem}\label{thm2.2} 
Fix $r$, $1\le r\le n$, and let $h=(h_1,\dots,h_r)$,
  where $0\le h_1<\dots<h_r\le n$ be given. Then,
\begin{equation}\label{2.8}
\mathbb{P}[p(Z_r(\tau))=h]=r!\mathbb{P}_{\text{Kr\,},N,K,q}[h],
\end{equation}
where $q=w^2(1+w^2)^{-1}$. Hence, if we disregard the ordering of the
particles in the zig-zag particle configuration, the probability of
$h$ is exactly $\mathbb{P}_{\text{Kr\,},N,K,q}[h]$.
\end{theorem}
\begin{Proof}
By lemma \ref{lem2.1} we have type I DR-paths from $(r+1-j,0)$ to
$(r,h_j)$, $1\le j\le r$, in CS-I, and type II DR-paths from $(j,0)$
to $(n+1-r,n-k_j)$, $1\le j\le n+1-r$ in CS-II, where
$k_1<\dots<k_{n+1-r}$ and (\ref{2.2}) holds. Together these describe
the whole domino tiling. Let $\omega[h]$ be the weight of all the type
I DR-paths between the specified points and $\omega'[h]$ the weight of
all the type II DR-paths between the given points. Then,
\begin{equation}\label{2.9}
\mathbb{P}[p(Z_r(\tau))=h]=\frac{\omega[h]\omega'[h]}{\sum_{0\le
    h_1<\dots<h_r\le n}\omega[h]\omega'[h]}.
\end{equation}
The quantities $\omega[h]$ and $\omega'[h]$ can be computed using the
Lindstr\"om-Gessel-Viennot method, \cite{Li}, \cite{GV}, see also 
\cite{Stem}, which
is a development of the Karlin-McGregor result in a discrete
setting. We want to compute the weight of a path that takes $n$ steps
to the right and $m$ steps up. Let $a$ be the number of $(1,0)$ steps,
$b$ the number of $(0,1)$ steps and $c$ the number of $(1,1)$
steps. Then, $n=a+c$ and $m=b+c$. The number of paths with a given
number of steps $a,b,c$ equals
\begin{equation}
\frac{(a+b+c)!}{a!b!c!}=\frac{(n+m-c)!}{(n-c)!(m-c)!c!}.\notag
\end{equation}
Note that $c$ can take all values between $0$ and $\min(n,m)$. The
total weight of all possible paths from $(0,0)$ to $(n,m)$ is thus
\begin{equation}
w(n,m)=\sum_{c=0}^{\min(n,m)}\frac{(n+m-c)!}{(n-c)!(m-c)!c!} w^{n+m-2c}.
\end{equation}
If we use the Pochhammer symbol $(a)_k=a(a+1)\dots(a+k-1)$, $(a)_0=1$,
then
\begin{equation}\label{2.10}
w(n,m)=\frac{w^{n+m}}{n!}\sum_{c=0}^\infty (m-c+1)_m(n-c+1)_c\frac
  {w^{-2c}}{c!}.
\end{equation}
The Lindstr\"om-Gessel-Viennot method now shows that the total wight
of all possible $r$ non-intersecting paths from $(r+1-j,0)$ to
$(r,h_j)$, $1\le j\le r$, is
\begin{equation}\label{2.11}
\omega[h]=\det(w(i-1,h_j))_{i,j=1}^r.
\end{equation}
Similarly, the total weight of all possible $n+1-r$ non-intersecting
paths from $(j,0)$ to $(n+1-r,n-k_j)$ is
\begin{align}
\omega'[h]&=\det(w(n+1-r-i,n-k_j))_{i,j=1}^{n+1-r}\notag\\
&=\det(w(i-1,k_{n+2-r-j}))_{i,j=1}^{n+1-r}.\notag
\end{align}
If we set $r'=n+1-r$, $h_j'=n-k_{n+2-r-j}$, then
\begin{equation}\label{2.12}
\omega'[h]=\det(w(i-1,h_j'))_{i,j=1}^{r'},
\end{equation}
which has exactly the same form as (\ref{2.11}). This is the advantage
of using both types of DR-paths.

Now, 
\begin{equation}\label{2.13}
\det(w(i-1,h_j))_{i,j=1}^r=\left(\prod_{j=1}^r
\frac{(1+w^2)^{j-1}}{w^{j-1}(j-1)!}\right)\Delta_r(h)\prod_{j=1}^r w^{x_j}.
\end{equation}
To see this, insert (\ref{2.10}) into the left hand side of
(\ref{2.13}), and use the multilinearity of the determinant to obtain
\begin{align}
&\left(\prod_{j=1}^{r}\frac{w^{j-1+h_j}}{(j-1)!}\right)
\sum_{c_1,\dots,c_r=0}^\infty
\prod_{i=1} (i-c_i)_{c_i}\frac
1{c_i!w^{2c_i}}\det((h_j-c_i+1)_{i-1})_{i,j=1}^r \notag\\
&=\left(\prod_{j=1}^{r}\frac{w^{j-1+h_j}}{(j-1)!}\right)\Delta_r(h)
\prod_{i=1}^r\sum_{c=0}^\infty (i-c)_c\frac 1{c! w^{2c}},\notag
\end{align}
which equals the right hand side of (\ref{2.13}) since,
\begin{equation}
\sum_{c=0}^\infty (i-c)_c\frac 1{c! w^{2c}}=\sum_{c=0}^{i-1}
\binom{i-1}{c}\frac 1{w^{2c}}=(1+\frac 1{w^2})^{i-1}.\notag
\end{equation}
From (\ref{2.12}) we obtain, after some manipulation,
\begin{equation}\label{2.14}
\omega'[h]=\left(\prod_{j=1}^{n+1-r}\frac{(1+w^2)^{j-1}}{w^{j-1}(j-1)!}\right) \Delta_{n+1-r}(k)\prod_{j=1}^{n+1-r}w^{n-k_j}.\notag
\end{equation}
Lemma 2.2 in \cite{Jo2} shows that if $s_1<\dots<s_N$ and
$r_1<\dots<r_M$ and the union of these two sets of numbers is exactly
$\{0,\dots,N+M-1\}$, then
\begin{equation}\label{2.13'}
\Delta_M(r)=\left(\prod_{j=1}^{N+M-1}j! 
\right)\left(\prod_{j=1}^N\frac 1{s_j!(N+M-1-s_j)!}\right)\Delta_N(s).
\end{equation}
If we use this formula, we obtain
\begin{equation}
\omega[h]\omega'[h]=w^{-r(r-1)}(1+w^2)^{n(n+1)/2-nr+r(r-1)}\prod_{j=1}^{r-1}
\frac{(n-j)!}{n!j!}\Delta_r(h)^2\prod_{j=1}^r\binom{n}{h_j} w^{2h_j}.\notag
\end{equation}
By this formula and (\ref{2.5}) we obtain
\begin{equation}
\omega[h]\omega'[h]=(1+w^2)^{n(n+1)/2}\frac {r!}{Z_{r,n,q}}\Delta_r(h)^2\prod_{j=1}^r\binom{n}{h_j}q^{h_j}p^{n-h_j},\notag
\end{equation}
where $q=w^2(1+w^2)^{-1}$. Hence
\begin{equation}\label{2.15}
\sum_{0\le h_1<\dots<h_r\le n}\omega[h]\omega'[h]=(1+w^2)^{n(n+1)/2},
\end{equation}
and the theorem follows from (\ref{2.9}).
\end{Proof}

As a corollary we obtain the following result first proved in
\cite{EKLP}, see also \cite{JPS}.

\begin{corollary}\label{corr2.3}
 The number of elements in $\mathcal{T}(A_n)$ is
  $2^{n(n+1)/2}$, and the probability of having $2k$ vertical tiles,
  $0\le k\le n(n+1)/2$, is
\begin{equation}\label{2.16}
\binom{n(n+1)/2}{k}\left(\frac {w^2}{1+w^2}\right)^{k}\left(\frac
  1{1+w^2}\right)^{n(n+1)/2-k}.
\end{equation}
The probability of having an odd number of vertical tiles is zero.
\end{corollary}
\begin{Proof} The first result follows by putting $w=1$ in
  (\ref{2.15}), and (\ref{2.16}) follows by expanding the right hand
  side of (\ref{2.15}) using the binomial theorem.
\end{Proof}

\subsection{Asymptotics in the Krawtchouk ensemble}
There is an {\it equilibrium measure} associated with the Krawtchouk
ensemble, see sect. 2.2. in \cite{Jo1}, and sect. 4.2 below for some
more details. If $r\to\infty$ and $n\to\infty$ in such a way that
$r/n\to t\in (0,1)$, then the expectation of the 
discrete measure $\frac 1r\sum_{j=1}^r
\delta_{h_j/n}$ converges weakly to the equilibrium measure
$u_{t,q}(x)dx$, i.e. the equilibrium measure gives the asymptotic
distribution of the particles. The equilibrium distribution is scaled
so that its support is contained in $[0,1]$. The equilibrium measure
also gives the asymptotic distribution of the zeroes of the Krawtchouk
polynomials scaled to $[0,1]$. See \cite{DS2} for this result and
explicit formulas for the equilibrium measure and its support. Since
the position of the rightmost particle determines the boundary of the
north polar zone (and we can make an analogous analysis for the other
polar zones or use symmetry), we can use the equilibrium measure to
prove the arctic ellipse theorem, \cite{JPS} and \cite{CEP}. The
arctic ellipse theorem of Jockush, Propp and Shor says that the
boundary of the temperate zone, scaled by $1/n$, converges in probability
to an ellipse as $n\to\infty$.
\begin{theorem} \label{thm2.4}
If we scale the Aztec diamond by $1/n$ in the original
  coordinate system, then the boundary $\partial T_n$ of the temperate
  zone in a rescaled andom tiling of $A_n$ under the probability measure
  (\ref{2.1}), converges in probability as $n\to\infty$ , $r/n\to t\in
  (0,1)$, to the ellipse $E$,
\begin{equation}\label{2.17}
\frac{x^2}p+\frac{y^2}q=1,
\end{equation}
in the sense that $\mathbb{P}[\text{dist\,}(\partial
T_n,E)\ge\epsilon]\to 0$ for any fixed $\epsilon>0$. Let
$\text{dist\,}_I(\partial T_n,E)$ be the maximal distance from a point
on $\partial T_n$ inside $E$ to $E$, and $\text{dist\,}_O(\partial
T_n,E)$ be the same thing but from a point outside $E$. Given
$\epsilon>0$, there are positive constants $I(\epsilon)$ and
$J(\epsilon)$ such that
\begin{equation}
\limsup_{n\to\infty}\frac 1{n^2}\log\mathbb{P}[\text{dist\,}_I(\partial
T_n,E)\ge\epsilon]\le -I(\epsilon)\notag
\end{equation}
and
\begin{equation}
\limsup_{n\to\infty}\frac 1{n}\log\mathbb{P}[\text{dist\,}_O(\partial
T_n,E)\ge\epsilon]\le -J(\epsilon)\notag
\end{equation}
\end{theorem}
\begin{Proof}
We will indicate how the shape (\ref{2.17}) is
obtained. The large deviation formulas, which imply the convergence in
probability to the ellipse, follow from theorem 2.2 in \cite{Jo1}. See
sect. 4.2  below for some more details in the analogous result
for rhombus tilings of a hexagon. Let $(x_I,y_I)$ be coordinates in CS-I
and $(x,y)$ coordinates in the original coordinate system. Then,
$x=n+1/2-x_I-y_I$, $y=1/2-x_I+y_I$. It follows from \cite{DS2} that the
right endpoint of the support of the equilibrium measure is given by
$\beta(t,q)=tp+(1-t)q+2\sqrt{t(1-t)pq}$ if $0<t\le 1-q$. The
connection between the position of the rightmost particle and the
boundary of the north polar zone described above imply that the limiting
boundary of the north polar zone must be the curve $(0,1-q]\ni t\to
(t,\beta(t,q))=(x_I,y_I)$ in CS-I. Using the coordinate transformation
we find that, in the original coordinate system, points on this curve
satisfy (\ref{2.17}).
\end{Proof}
Let $\nu[a,b](h)$ be the number of particles in the interval $[a,b]$
in the particle configuration $h$. From the formula (\ref{2.3}) and
$h(Q_{n+1}^r)=2r-1$, we obtain
\begin{align}
&h(Q_k^r)-(2r-1)=\sum_{j=k}^n [h(Q_j^r)-h(Q_{j+1}^r)]\notag \\
&=-2\nu[0,n-k]+2(n-k+1-\nu[0,n-k])=2(n-k+1)-4\nu[0,n-k].\notag
\end{align}
Consequently,
\begin{equation}\label{2.18}
h(Q_k^r)=2(n-k+r)+1-4\nu[0,n-k].
\end{equation}
Assume that $k/n\to\tau$ and $r/n\to t$, $0<\tau<t$ as
$n\to\infty$. Then, by (2.18) and the weak convergence of the particle
distribution, 
\begin{equation}\label{2.19}
\lim_{n\to\infty}\frac{h(Q_k^r)}{n}=2(1-\tau+t)-4t\int_0^{1-\tau}
u_{t,q}(x)dx. 
\end{equation}
The precise form of the equilibrium measure is given in \cite{DS2} and
using this we can work out the asymptotic height function. We will not
evaluate these integrals here. We can also obtain large deviation
formulas (and estimates, compare lemma 4.1 in \cite{Jo1}) for
macroscopic deviations from the asymptotic (average) height
function. This is analogous to the large deviation formulas for the
Wigner semi-circle law, \cite{BG}. We will not develop the details,
since it is very analogous to the corresponding random matrix
results. See \cite{CEP} for previous large deviation estimates 
and asymptotics for the
height function.

We will now analyze the fluctuations of the height function, or what
is equivalent, by (\ref{2.18}), the fluctuations in the number of
particles in an interval in the Krawtchouk ensemble (\ref{2.4}). Let
$I=\{b-L,b-L+1,\dots,b\}\subseteq\{0,\dots,K\}$ be an ``interval'' of
  length $L$ and let $\nu(I)$ denote the number of particles in $I$,
\begin{equation}
\nu(I)=\#\{h_i\,;\,1\le i\le N\,\,\text{and}\,\, h_i\in I\},\notag
\end{equation}
We want to prove the following result for the variance of $\nu(I)$,
the {\it number variance}. 
\begin{proposition}\label{prop2.5}
Assume that $N/K\to t\in(0,1/2]$, $p=q=1/2$, $b/K\to\beta$,
$(b-L)/K\to\beta'\le\beta$ as $N,K,L\to\infty$, where $\beta$ or
$\beta'$ belongs to the interior, $S$, of the support of the
equilibrium measure $u_{t,1/2}$. Then
\begin{equation}\label{2.20}
\lim\frac{\text{var\,}(\nu(I))}{\log L}=\frac 1{\pi^2}\xi_{\beta,\beta'},
\end{equation}
where $\xi_{\beta,\beta'}$ is $=1$ if both $\beta$ and $\beta'$ belong
to S and $=1/2$ otherwise.
\end{proposition}

This type of result was conjectured in \cite{PrSt}.
We will give the proof of this proposition below, which is rather
long. 
The case $p\neq q$ could also be worked out, but we stick with
$p=q=1/2$ for simplicity.
Once we have this result we can apply the
Costin/Lebowtiz/Soshnikov argument, \cite{CL}, \cite{So1}, to prove
that the fluctuations are normal. This type of results have been
proved in other tiling models by Kenyon, see \cite{Ke1}, \cite{Ke2},
\cite{Ke3}.

\begin{theorem}\label{thm2.6} With the same assumptions as in the
proposition,
\begin{equation}\label{2.21}
\frac{\nu(I)-\mathbb{E}[\nu(I)]}{\sqrt{\text{var\,}(I)}}\Rightarrow N(0,1)
\end{equation}
as $K\to\infty$. i.e. we have convergence in distribution to a
standard normal random variable.
\end{theorem}
\begin{Proof}By a theorem in \cite{So1}, p. 8, see also \cite{So2}, the
  result (\ref{2.21}) follows from (\ref{2.20}) since we have a
  determinantal random point field, i.e. the correlation functions are
  given by determinants as in (\ref{2.6'}). The kernel (\ref{2.7})
  defines a trace class operator (it has finite rank) $\mathcal{K}$,
  which satisfies $0\le\mathcal{K}\le I$. This follows from
  (\ref{2.7}) and the orthogonality $(\ref{2.5'})$. Hence the
  conditions in the Costin-Lebowitz-Soshnikov theorem are satisfied.
\end{Proof}

Combining this theorem with (\ref{2.18}) we obtain the next theorem.

\begin{theorem}\label{thm2.7}
Take the uniform distribution on $\mathcal{T}(A_n)$ and $0\le r\le
n/2$ (the case $n/2\le r\le n$ is similar by symmetry). Let
$Q_k^r=(-r+k, n+1-k-r)$ as before and let $h(Q_k^r)$ be the value of
the domino height function above this point. If $r/n\to t\in (0,1)$,
$k/n\to\kappa \in(1/2-\sqrt{t(1-t)},1/2+\sqrt{t(1-t)})\doteq U_t$ or
$j/n\to\kappa'\in U_t$, and $|k-j|\to\infty$, then
\begin{equation}\label{2.22}
\frac{h(Q_k^r)-h(Q_j^r)-\mathbb{E}[h(Q_k^r)-h(Q_j^r)]}
{4\sqrt{\xi_{\kappa,\kappa'}\pi^{-2}\log|k-j|}}\Rightarrow N(0,1).
\end{equation}
Here $\xi_{\kappa,\kappa'}=1$ if $\kappa,\kappa'\in U_t$ and 
$\xi_{\kappa,\kappa'}=1/2$ if one of $\kappa$ or $\kappa'$ does not
belong to $U_t$.
\end{theorem}
Note that the condition on $\kappa,\kappa'$ corresponds exactly to the
condition that one of $\beta$ or $\beta'$ in proposition 2.5 belongs
to the interior of the support of the equilibrium measure.

We turn now to the proof of proposition 2.5 which is rather
lengthy. The proof is based on sufficiently good asymptotic control of
the Krawtchouk kernel, (\ref{2.7}). 
For results on asymptotics of Krawtchouk polynomials see \cite{IS}.
We state the needed result as a
lemma, which we will prove later. 
Let $\rho(\xi)$ be defined by (\ref{xx}) below. We consider the case when 
$p=1/2$ and $N/K\to t$. Then, $\rho'(\xi)=u_{t,1/2}(\xi)$ is the
equilibrium measure, which is suppported in 
$[1/2-\sqrt{t(1-t)},1/2+\sqrt{t(1-t)}]$, see \cite{DS2}. A computation gives,
\begin{equation}\label{2.22'}
\rho'(\xi)=\frac
1{\pi}\arctan\frac{\sqrt{t(1-t)-(\xi-1/2)^2}}{\sqrt{1/4-t(1-t)}}.
\end{equation}

\begin{lemma}\label{lem2.8}
Consider the Krawtchouk kernel (\ref{2.7}) with $p=1/2$ and let
$\delta>0$ (small). Set $t=N/K$ and assume $0<t\le 1/2$. If
$|x/K-1/2|\le \sqrt{t(1-t)}-\delta$ and $|y/K-1/2|\le
\sqrt{t(1-t)}-\delta$, then there is a constant $C$, independent of
$N,K,x$ and $y$ such that
\begin{equation}\label{2.23}
\left|(x-y)K_{\text{Kr\,},N,K,1/2}(x,y)-
\frac{\sin\pi K(\rho(x/K)-\rho(y/K))}{\pi}\right|\le C\left(\frac
1{\sqrt{K}}+ \frac{|x-y|}K\right).
\end{equation}
Also, if $|x/K-1/2|\le \sqrt{t(1-t)}-\delta$, there is a constant $C$
such that for $y\ge x$,
\begin{equation}\label{2.24}
|(x-y)K_{\text{Kr\,},N,K,1/2}(x,y)|\le C\min(K^{1/4},
 |t(1-t)-(y/K-1/2)^2|^{-1/2}) \frac{K^{1/4}}{(K-y)^{1/4}}.
\end{equation}
There is also an analogous result for $y\le x$. If
$1/2+\sqrt{t(1-t)}<1$, we can obtain a much better estimate for $y$
outside the support of the equilibrium measure; compare (\ref{2.59}) below.
\end{lemma}

We will now prove proposition 2.5.
\begin{Proof}
Write $\mathcal{K}(x,y)=K_{\text{Kr\,},N,K,1/2}(x,y)$. It follows from
(\ref{2.7}) and the orthogonality (\ref{2.5'}), that
$\mathcal{K}(x,y)$ is a reproducing kernel
\begin{equation}\label{2.25}
\sum_{j=0}^K\mathcal{K}(i,j)\mathcal{K}(j,k)=\mathcal{K}(i,k).
\end{equation}
Set $\tilde{I}=\{0,\dots,K\}\setminus I$. Then by the formulas
(\ref{2.6'}) and (\ref{2.25}), we see that
\begin{align}\label{2.26}
\text{var\,}[\nu(I)]&=\sum_{j\in I}\mathcal{K}(j,j)-\sum_{i,j\in I}
\mathcal{K}(i,j)^2 \notag \\
&=\sum_{j\in
  I}\left(\sum_{i=0}^K\mathcal{K}(j,i)\mathcal{K}(i,j)\right)
-\sum_{i,j\in I} \mathcal{K}(i,j)^2 \notag \\
&=\sum_{j\in
  I}\sum_{i\in\tilde{I}}\mathcal{K}(i,j)^2=\Sigma_1+\Sigma_2,
\end{align}
where
\begin{align}
\Sigma_1&=\sum_{i=0}^L\sum_{j=1}^{K-b}\mathcal{K}(b-i,b+j)^2\notag\\
\Sigma_1&=\sum_{i=0}^{b-L}\sum_{j=1}^{L}\mathcal{K}(b-L-i,b-L+j)^2.\notag
\end{align}
We will consider the case when $\beta$ lies in the support $S$ of the
equilibrium measure. The contribution to $\Sigma_1$ comes from the
right endpoint of the interval, and the 
contribution to $\Sigma_2$ from the left endpoint. We will show that
\begin{equation}\label{2.27}
\Sigma_1=\frac 1{2\pi^2}\log L+ O(\log(\log L)).
\end{equation}
The same thing is true for $\Sigma_2$, with an analogous proof, if
$\beta'\in S$, so if we establish (\ref{2.27}) the proposition is
proved.

Let $a(L)=[L/\log L]$. We split $\Sigma_1$ into two parts
\begin{align}\label{2.28}
\Sigma_1&=\sum_{i=0}^{L}\sum_{j=1}^{a(L)}\mathcal{K}(b-i,b+j)^2+
\sum_{i=0}^{L}\sum_{j=a(L)+1}^{K-b}\mathcal{K}(b-i,b+j)^2\notag\\
&=\Sigma_1'+\Sigma_1''.
\end{align}
Now, by (\ref{2.24}),
\begin{align}
\Sigma_1''&\le\sum_{i=0}^L\sum_{j=a(L)+1}^{[(\alpha-\delta/2)K]-b}
  \frac{CK^{1/2}}{|K\alpha-(b+j)|^{1/2}}\frac 1{(i+j)^2}\notag\\
&+\sum_{i=0}^\infty\sum_{j=[(\alpha-\delta/2)K]-b}^{[\alpha K]-[K^{1/2}]-b}
\frac{CK^{1/2}}{|K\alpha-(b+j)|^{1/2}}\frac 1{(i+j)^2}\notag\\
&+\sum_{i=0}^\infty\sum_{j=[\alpha K]-[K^{1/2}]-b}^{[\alpha
  K]-[K^{1/2}]+b}
\frac{CK^{1/4}}{(i+j)^2}\notag\\
&+\sum_{i=0}^\infty\sum_{j=[\alpha K]-[K^{1/2}]+b}^{K-b}
\frac{CK^{1/2}}{|K\alpha-(b+j)|^{1/2}}\frac 1{(i+j)^2}\notag,
\end{align}
where $\alpha=1/2+\sqrt{t(1-t)}$, the right endpoint of the
support. The first sum is
\begin{equation}
\le C\sum_{i=0}^L\sum_{j=a(L)+1}^\infty\frac 1{(i+j)^2}\le C\log (\log
L).\notag
\end{equation}
In the last three sums, we use $\sum_{i=0}^\infty 1/(i+j)^2\le
1/(j-1)$, and it is then easy to see that the $j$-sums are $\le
C$. Thus
\begin{equation}\label{2.29}
\Sigma_1''\le C(1+\log (\log L)).
\end{equation}
We also split $\Sigma_1'$ into two sums
\begin{align}\label{2.30}
\Sigma_1'&=\sum_{i=0}^{a(L)}\sum_{j=1}^{a(L)}\mathcal{K}(b-i,b+j)^2+
\sum_{i=a(L)+1}^{L}\sum_{j=1}^{a(L)}\mathcal{K}(b-i,b+j)^2\notag\\
&=S_1+S_2.
\end{align}
The second sum is estimated in the same way as $\Sigma_1''$ using
(\ref{2.24}), and this gives
\begin{equation}\label{2.31}
S_2\le C.
\end{equation}
To control $S_1$ we use (\ref{2.23}), which gives
\begin{align}\label{2.32}
&\left|S_1-\sum_{i=0}^{a(L)}\sum_{j=1}^{a(L)}\frac {\sin^2\pi
    K(\rho((b-i)/K)- \rho((b+j)/K))}{\pi^2(i+j)^2}\right|\notag\\
&\le C\sum_{i=0}^{a(L)}\sum_{j=1}^{a(L)}(\frac
    1{\sqrt{K}}+\frac{i+j}{K})\frac 1{(i+j)^2}\le C.
\end{align}
If we use $\sin^2x=(1+\cos 2x)/2$ and observe that
\begin{equation}\label{2.33}
\sum_{i=0}^{a(L)}\sum_{j=1}^{a(L)}\frac 1{(i+j)^2}=\log L+ O(\log
(\log L)),
\end{equation}
we see that what remains to be proved is
\begin{equation}\label{2.34}
\left|\sum_{i=0}^{a(L)}\sum_{j=1}^{a(L)}\frac {\cos 2\pi
    K(\rho((b-i)/K)- \rho((b+j)/K))}{\pi^2(i+j)^2}\right|\le C\log
    (\log L).
\end{equation}
When we have this we note that (\ref{2.28}) - (\ref{2.34}) imply
(\ref{2.27}) and we are done. 

To prove (\ref{2.34}) we use summation by parts. Set $S(0)=0$ and
\begin{equation}
S(j)=\sum_{n=1}^j \exp (-2\pi Ki\rho(\frac{b+n}K)),\quad 1\le j\le
a(L).\notag 
\end{equation}
Consider the expression
\begin{equation}\label{2.35}
\sum_{m=0}^{a(L)} e^{2\pi Ki\rho(\frac{b-m}K)}\sum_{j=1}^{a(L)}
\frac 1{(m+j)^2}e^{-2\pi Ki\rho(\frac{b+j}K)},
\end{equation}
whose real part is what we want.
Summation by parts gives
\begin{equation}\label{2.36}
\sum_{j=1}^{a(L)}
\frac 1{(m+j)^2}e^{-2\pi Ki\rho(\frac{b+j}K)}=
\frac 1{(m+a(L))^2}S(a(L))+\sum_{j=1}^{a(L)-1}\frac
{2(m+j)+1}{(m+j)^2(m+j+1)^2} S(j).
\end{equation}
Since $|S(a(L)|\le a(L)$, the first term in (\ref{2.36}) gives a
contribution $\le 1$ to the expression (\ref{2.35}). The contribution
of the second term in (\ref{2.36}) to (\ref{2.35}) is
\begin{equation}\label{2.37}
\le \sum_{m=0}^{a(L)}\sum_{j=1}^{a(L)-1}\frac
{2(m+j)+1}{(m+j)^2(m+j+1)^2} S(j).
\end{equation}
Fix an integer $\Delta>1$ and assume that $\Delta< j\le a(L)$. Write
$j=k\Delta+r$, where $0\le r<\Delta$, and $k\le j/\Delta$. Then
\begin{equation}\label{2.38}
S(j)=\sum_{u=1}^k\sum_{v=1}^\Delta e^{-2\pi
  Ki\rho(\frac{b+u\Delta+v}K)}+ \sum_{v=1}^re^{-2\pi
  Ki\rho(\frac{b+k\Delta+v}K)}.
\end{equation}
Write $\xi=(b+u\Delta)/K$. Then
$\rho(\frac{b+u\Delta+v}K)=\rho(\xi)+\frac
vK\rho'(\xi)+O(\frac{v^3}{K^2})$. Note that $\rho'(\xi)>0$ since we
are inside the support of the equilibrium measure $\rho'$. Also,
$\rho'(\xi)<1$ by (\ref{2.22'}). Using this we see that
\begin{equation}\label{2.39}
\sum_{v=1}^\Delta\left|e^{-2\pi Ki\rho(\xi+\frac vK)}-
e^{-2\pi Ki\rho(\xi)-2\pi iv\rho'(\xi)}\right|
\le C\sum_{v=1}^\Delta\frac{v^2}K\le C\frac{\Delta^3}K.
\end{equation}
Now,
\begin{equation}\label{2.40}
\left|\sum_{v=1}^\Delta e^{-2\pi Ki\rho(\xi)-2\pi iv\rho'(\xi)}\right|\le C,
\end{equation}
since $0<\rho'(\xi)<1$. Similarly
\begin{equation}\label{2.41}
\left|\sum_{v=1}^\Delta e^{-2\pi Ki\rho(\frac{b+k\Delta+v}K)}
\right|\le C(1+\frac{\Delta^3}K).
\end{equation}
Combining (\ref{2.38}) to (\ref{2.41}) we obtain
\begin{equation}
|S(j)|\le Ck(1+\frac{\Delta^3}K)\le Cj(\frac 1{\Delta}+\frac{\Delta^2}K)\notag
\end{equation}
for $\Delta <j\le a(L)$. Since evidently $|S(j)|\le j$ for $1\le j\le
\Delta$, we obtain the estimate
\begin{equation}
\left|\sum_{m=0}^{a(L)}\sum_{j=1}^{a(L)-1}\frac
{2(m+j)+1}{(m+j)^2(m+j+1)^2} S(j)\right|\le C\log\Delta+C(\frac
1{\Delta}+\frac{\Delta^2}K) \log L\notag
\end{equation}
of the expression in (\ref{2.37}). If we take $\Delta=[\log L]$, we
see that the absolute value of the expression (\ref{2.35}) is $\le
\log (\log L)$, which gives the bound (\ref{2.34}).
\end{Proof}
We still have to prove lemma 2.8.

\begin{Proof}
Let $p=q=1/2$ in (\ref{2.7}) and make the change of variable
$z=2w$. This gives,
\begin{equation}
p_n(z)=\binom{K}{n}^{-1/2}\frac 1{2\pi
  i}\int_{\gamma}\frac{(1+w)^x(1-w)^{K-x}}{w^n}\frac{dw}w.\notag
\end{equation}
Set
\begin{equation}\label{2.40'}
A_N(x,y)=\frac{\kappa_{N-1}}{\kappa_N}(w(x)w(y))^{1/2}\binom{K}{N}^{-1/2}
\binom{K}{N-1}^{-1/2},
\end{equation}
$G_N(z;x)=(1+z)^x(1-z)^{K-x}w^{-N}$ and
$L_N(x,y)=(x-y)K_{\text{Kr\,}, N, K, 1/2}(x,y)$. Then,
\begin{equation}
L_N(x,y)=A_N(x,y)\int_{\gamma}\frac{dz}{2\pi i
  z}\int_{\gamma}\frac{dw}{2\pi i w} G_N(z;x)G_N(w;y)(w-z).\notag
\end{equation}
Set $\xi=x/K$ and $\eta=y/K$ and define $f(z)=\xi\log (1+z)+(1-\xi
)\log(1-z)-t\log z$, which is the relevant function in the saddle
point argument. The equation $f'(z)=0$ has the solutions $z_c^\pm=
(1-t)^{-1}(\xi-1/2\pm i\sqrt{t(1-t)-(\xi-1/2)^2})$. Let
$r_c=|z_c^\pm|=\sqrt{t(1-t)^{-1}}$. We choose $\gamma$ to be the
circle with radius $r_c$, and write $z_c^\pm=r_c\exp(\pm
i\theta_c(\xi))$, where
\begin{equation}\label{2.41'}
\cos\theta_c(\xi)=\frac{\xi-1/2}{\sqrt{t(1-t)}},\quad
0\le\theta_c\le\pi.
\end{equation}
We assume that $\xi\ge 1/2$ (the other case is similar by
symmetry). Also, we will first assume that $\xi\le 1/2+\sqrt{t(1-t)}$,
i.e. we are inside the support of the equilibrium measure. It follows
that $0\le \theta_c\le \pi/2$.
We obtain,
\begin{equation}\label{2.42}
L_N(x,y)=A_N(x,y)\frac{r_c}{(2\pi)^2}\int_{-\pi}^\pi\int_{-\pi}^\pi
G_N(r_ce^{i\theta};x)G_N(r_ce^{i\phi};y)
(e^{i\phi}-e^{i\theta})d\theta d\phi,
\end{equation}
and this is the formula we will use. Let
\begin{align}
g(\theta)&=g(\theta;x)=\log|G_N(r_ce^{i\theta};x)|\notag\\
&=-N\log r_c+\frac x2\log|1+r_ce^{i\theta}|^2+
\frac{K-x}2\log|1-r_ce^{i\theta}|^2.\notag
\end{align}
Taking the derivative gives
\begin{equation}\label{2.43}
g'(\theta)=\frac{2Kt}{1-t}\sin\theta
\left[\frac{\cos\theta-\frac{\xi-1/2}{\sqrt{t(1-t)}}}
  {|1+r_ce^{i\theta}|^2 |1-r_ce^{i\theta}|^2}\right].
\end{equation}
From this formula we see that
\begin{equation}\label{2.44}
g(\theta)\le g(\theta_c),\quad -\pi\le\theta\le\pi.
\end{equation}
Taking the absolute values in (\ref{2.42}) gives
\begin{align}\label{2.45}
|L_N(x,y)|&\le\frac{2r_c}{\pi^2}A_N(x,y)e^{g(\theta_c(\xi);x)+
g(\theta_c(\eta);y)}\notag\\
&\times\left(\int_0^\pi
  e^{g(\theta;x)-g(\theta_c(\xi);x)}d\theta\right)
\left(\int_0^\pi e^{g(\theta;y)-g(\theta_c(\eta);y)}d\theta\right).
\end{align}
Fix $\delta >0$ (small). If $\delta\le\theta_c\le\pi/2$ we can make a
quadratic approximation around $\theta_c$ to obtain the estimate
\begin{equation}\label{2.46}
\int_0^\pi e^{g(\theta)-g(\theta_c)}d\theta\le\frac{C}{\sqrt{K}},
\end{equation}
for some constant $C$. If $\theta_c$ is close to $0$ we have to be
more careful. A computation using (\ref{2.43}) shows that there is a
constant $\alpha>0$ such that
\begin{equation}\label{2.47}
\int_0^\pi e^{g(\theta)-g(\theta_c)}d\theta
\le\int_0^\infty e^{-\alpha K(u^2-u_c^2)^2}du.
\end{equation}
Here we have taken $u=\sin\theta$, $u_c=\sin\theta_c$. Now,
\begin{equation}\label{2.48}
\int_0^\infty e^{-\alpha K(u^2-u_c^2)^2}du\le C\min(\frac
1{u_c\sqrt{K}},\frac 1{K^{1/4}})
\end{equation}
and hence
\begin{equation}\label{2.49}
\int_0^\pi e^{g(\theta;x)-g(\theta_c(\xi);x)}d\theta\le
C\min(\frac
1{\sin\theta_c(\xi)\sqrt{K}},\frac 1{K^{1/4}}).
\end{equation}
To prove (\ref{2.48}) assume that $u_c\ge 1/K^{1/4}$. Set
$u=u_c(1+s)$. The left hand side becomes
\begin{equation}
u_c\int_{-1}^\infty e^{-\alpha Ku_c^4s^2(s+2)^2}ds\le
u_c\int_{-\infty}^\infty e^{-\alpha Ku_c^4s^2}ds=\frac C{u_c\sqrt{K}}.
\notag
\end{equation}
On the other hand if $u_c\le 1/K^{1/4}$, then the left hand side
equals
\begin{equation}
\int_0^{u_c\sqrt{2}}e^{-\alpha K(u^2-u_c^2)^2}du+
\int_{u_c\sqrt{2}}^\infty e^{-\alpha K(u^2-u_c^2)^2}d\le\sqrt{2}u_c
+\int_0^\infty e^{-\alpha Ku^4}du\le\frac{C}{K^{1/4}},
\notag
\end{equation}
and we have proved (\ref{2.48}). Note that $\sin\theta_c(\xi)=\frac
1{\sqrt{t(1-t)}}\sqrt{t(1-t)-(\xi-1/2)^2}$, so (\ref{2.49}) can be
written
\begin{equation}\label{2.50}
\int_0^\pi e^{g(\theta;x)-g(\theta_c(\xi);x)}d\theta\le C\min(\frac
1{\sqrt{K} \sqrt{|t(1-t)-(\xi-1/2)^2|}},\frac 1{K^{1/4}}).
\end{equation}
If $\xi\ge1/2+\sqrt{t(1-t)}$ we can take $\theta_c=0$ and then
(\ref{2.50}) still holds by a similar argument using
(\ref{2.43}). Thus, for any $0\le x,y\le K$, by (\ref{2.45}) and
(\ref{2.50}),
\begin{align}\label{2.51}
|L_N(x,y)|&\le CA_N^\ast(x,y)\min(|t(1-t)-(\xi-1/2)^2|^{-1/2},
K^{1/4})\notag\\ 
&\times \min(|t(1-t)-(\eta-1/2)^2|^{-1/2}, K^{1/4}),
\end{align}
where $\theta_c(\xi)$ is given by (\ref{2.41'}) if $|\xi-1/2|\le
\sqrt{t(1-t)}$, $\theta_c(\xi)=0$ if $|\xi-1/2|>\sqrt{t(1-t)}$ and
\begin{equation}
A_N^\ast(x,y)=\frac 1Kr_c
A_N(x,y)e^{g(\theta_c(\xi);x)+g(\theta_c(\eta);y)}.
\notag
\end{equation}

Let $\delta>0$ and consider $x,y$ such that $|\xi-1/2|\le
\sqrt{t(1-t)}-\delta$ and $|\eta-1/2|\le
\sqrt{t(1-t)}-\delta$, i.e. we are inside the support of the
equilibrium measure. Now, 
\begin{equation}
\left.\frac{d^2}{d\theta^2}f(r_ce^{i\theta})\right|_{\theta=\theta_c}
=-(z_c^+)\doteq -d_+(\xi)\notag
\end{equation}
and a straightforward computation shows that $\Re
d_+(\xi)>0$. Furthermore, if $d_-(\xi)\doteq 
-\frac{d^2}{d\theta^2}f(r_ce^{i\theta})|_{\theta=-\theta_c}$, then
$d_-(\xi)=\overline{d_+(\xi)}$. Also, for the $\xi$ and $\eta$ we are
considering we have that $|d_+(\xi)-d_+(\eta)|\le C|\xi-\eta|$. A
standard  local saddle-point argument now gives, using (\ref{2.42}),
(\ref{2.43}) and (\ref{2.44}),
\begin{equation}
L_N(x,y)=\frac {r_cA_N(x,y)}{2\pi}(I_{++}+I_{-+}+I_{+-}+I_{--}),\notag
\end{equation}
where ($a,b=\pm$),
\begin{equation}
I_{ab}=\left[\frac 1{K\sqrt{d_a(\xi)d_b(\eta)}}+O(\frac
  1{K^{3/2}})\right] G_N(r_ce^{ai\theta_c(\xi)};x)
G_N(r_ce^{ai\theta_c(\eta)};y)(e^{bi\theta_c(\eta)}-e^{ai\theta_c(\xi)}).
\notag
\end{equation}
Write,
\begin{equation}\label{xx}
G_N(r_ce^{ai\theta_c(\xi)};x)=e^{g(\theta_c(\xi);x)+iK\pi\rho(\xi)},
\end{equation}
which defines $\rho$ in lemma \ref{lem2.8}.
Since $|\exp(i\theta_c(\eta))-\exp(i\theta_c(\xi))|\le C|\xi-\eta|$, 
$d_+(\xi)=\overline{d_-(\xi)}$ and $|d_+(\xi)-d_+(\eta)|\le
C|\xi-\eta|$, we get
\begin{align}
&\left|L_N(x,y)-\frac{A_N^\ast(x,y)}{2\pi|d_+(\xi)|}\left[
e^{iK(\rho(\xi)-\rho(\eta))}(e^{-i\theta_c(\xi)}-e^{i\theta_c(\eta)})\right.
\right.\notag\\&+\left.\left.
e^{iK(-\rho(\xi)+\rho(\eta))}(e^{i\theta_c(\xi)}-e^{-i\theta_c(\eta)})\right] 
\right|
\le CA_N^\ast(x,y)(\frac 1{\sqrt{K}}+|\xi-\eta|),
\notag
\end{align}
which can be written
\begin{equation}\label{2.52}
\left|L_N(x,y)-\frac{2\sin\theta_c(\xi)}{\pi|d_+(\xi)|}A_N^\ast(x,y)
\sin\pi K(\rho(\xi)-\rho(\eta))\right|\le 
CA_N^\ast(x,y)(\frac 1{\sqrt{K}}+\frac{|x-y|}K).
\end{equation}

We now investigate $A_N^\ast(x,y)$ and start with the case when both
$\xi$ and $\eta$ are inside the support of the equilibrium
measure. Inserting the formulas for $\kappa_n$ and $w(x)$ in
(\ref{2.40'}) we obtain
\begin{equation}\label{2.53}
A_N(x,y)=(K-N+1)\binom{K}{N}^{-1}2^{-K}\binom{K}{x}^{1/2}\binom{K}{y}^{1/2}.
\end{equation}
Furthermore, a computation shows that
\begin{align}\label{2.54}
e^{g(\theta_c(\xi);x)}&=r_c^{-N}(1+r_c^2+2r_c\cos\theta_c)^{x/2}
(1+r_c^2-2r_c\cos\theta_c)^{(K-x)/2}\notag\\
&=\frac{(1-t)^{N/2}2^{K/2}}{t^{N/2}(1-t)^{K/2}K^{K/2}}x^{x/2}(K-x)^{(K-x)/2}.
\end{align}
Stirling's formula gives the asymptotic formulas
\begin{equation}\label{2.55}
\binom{K}{x}=\frac{K^K}{(K-x)^{K-x}x^x\sqrt{K}}\frac 1{\sqrt{2\pi\xi
    (1-\xi)}} (1+O(\frac 1{K}))
\end{equation}
and
\begin{equation}\label{2.56}
\binom{K}{N}=\frac 1{(1-t)^{K-N}t^N\sqrt{K}}\frac 1{\sqrt{2\pi t
    (1-t)}}(1+O(\frac 1{K})).
\end{equation}
Combining (\ref{2.53}) - (\ref{2.56}) gives
\begin{equation}\label{2.57}
A_N^\ast(x,y)=\frac{t(1-t)}{2\sqrt{\xi(1-\xi)}}(1+O(\frac 1K)+O(|\xi-\eta|)).
\end{equation}
A computation shows that
\begin{equation}\label{2.58}
\frac{\sin\theta_c}{|d_+(\xi)|}=\frac{\sqrt{\xi(1-\xi)}}{t(1-t)},
\end{equation}
and hence (\ref{2.23}) in the proposition follows by combining
(\ref{2.52}), (\ref{2.57}) and (\ref{2.58}).

We also have to estimate $A_N^\ast(x,y)$ when $\xi$ or $\eta$ is
outside the support of the equilibrium measure. Assume that $\xi$ is
inside and $\eta$ outside the support (or close to the edge of the
support), so that $\theta_c(\eta)=0$. We have that
\begin{equation}
e^{g(0;y)}=r_c^{-N}(1+r_c)^y(1-r_c)^{K-y}.\notag
\end{equation}
Set
\begin{equation}
B(y)=r_c^{-N}\left(1+r_c^2+2r_c\frac{\eta-1/2}{\sqrt{t(1-t)}}\right)^{y/2}
\left(1+r_c^2-2r_c\frac{\eta-1/2}{\sqrt{t(1-t)}}\right)^{(K-y)/2},
\notag
\end{equation}
which corresponds to the expression (\ref{2.54}) with $x=y$. We want
to show that
\begin{equation}\label{2.59}
\frac{e^{g(0;y)}}{B(y)}\le\gamma^K\le 1,
\end{equation}
with $\gamma<1$ if $\eta>1/2+\sqrt{t(1-t)}$. A computation gives
\begin{equation}
\frac{e^{g(0;y)}}{B(y)}=\left(\frac{\alpha}{\eta}\right)^{y/2}
\left(\frac{1-\alpha}{1-\eta}\right)^{(K-y)/2}=
\left[e^{\eta\log\frac{\alpha}{\eta}+(1-\eta)\log(\frac{1-\alpha}{1-\eta})} 
\right]^K,\notag
\end{equation}
where $\alpha=1/2+\sqrt{t(1-t)}$, $0\le\alpha\le\eta\le 1$. Set
$\gamma=e^{\eta\log\frac{\alpha}{\eta}+
(1-\eta)\log(\frac{1-\alpha}{1-\eta})}$ and note that
$\eta\log\frac{\alpha}{\eta}+(1-\eta)\log(\frac{1-\alpha}{1-\eta})\le
\log(\alpha+1-\alpha)=0$ by convexity; if $\eta>\alpha$ we get a
strict inequality. When $y$ is close to $K$ we have to be somewhat more
careful in estimating the binomial coefficients. We use
\begin{equation}
\binom{K}{y}=\frac{K^K}{(K-y)^{K-y}y^y}\sqrt{\frac{K}{y(K-y)}}(1+O(\frac
1y) +O(\frac 1{K-y})).
\notag
\end{equation}
Using this formula and proceeding as before we obtain
\begin{equation}\label{2.60}
A_N^\ast(x,y)=\frac{t(1-t)}{(\xi(1-\xi))^{1/4}}\frac{K^{1/4}}{(K-y)^{1/4}}
(1+O(\frac 1{|K-y|})+O(\frac 1K)),
\end{equation}
instead of (\ref{2.57}). Combining (\ref{2.51}) and (\ref{2.60}) we
obtain (\ref{2.24}) and the proposition is proved. Note that when
$t=1/2$ some modifications in the arguments above are needed. We will
omit the details.
\end{Proof}
\subsection{The corner growth model}

We can draw the type I DR-paths in a different way so that they look
like the heights of a cascade of discrete polynuclear growth (PNG)
models. If we place a W-domino so that it has corners at $(0,0),
(1,0), (1,2)$ and $(0,2)$ we draw a path by connecting the points
$(0,1/2)$, $(1/2,1/2)$, $(1/2,3/2)$ and $(1,3/2)$ with straight line
segments. We draw a path on an E-domino analogously, it is the mirror
image of the W-domino in the middle vertical line. In this way, from
the DR-paths,
we obtain height curves, which we can think of as graphs of functions
$h_k(x,n)$, $1\le k\le n$, where the $k\,$:th curve goes from
$E_k=(-n-1+k,1/2-k)$ to $A_k=(n+1-k,1/2-k)$ in the original coordinate
system. We think of $h_k(x,n)$ as the height of the level-$k$ growth
process at the point $x$ at time $t=n$. The level-$k$ height curve
does not intersect the level-$(k+1)$ height curve. Note that the
vertical steps always have step size $\pm 1$;
$h_k(x+,n)-h_k(x-,n)\in\{-1,0,1\}$, and the jumps can occur only at
points $x\in\{-n-1/2+k,-n+1/2+k,\dots,n+1/2-k\}$.

A random tiling of the Aztec diamond can be generated by the so called
{\it shuffling algorithm} which we will now describe briefly, see
\cite{EKLP} and \cite{JPS} for more details. Start with $A_1$. We can
tile it by either two horizontal dominoes, with probability $1-q$, or
two vertical dominoes, with probability $q$, (compare
(\ref{2.16}). Assume that we have a random tiling of $A_k$ for some
$k\ge 1$. We will define a random tiling of $A_{k+1}$, given a tiling
of $A_k$. Call two dominoes which share a side of length two a {\it
  pair}. Two horizontal dominoes form a bad pair if the lower one is
N and the upper one is S. Similarly, a pair of vertical dominoes is
bad if the left one is E and the right one is W. All other pairs are
{\it good}. In the first step we remove all bad pairs of dominoes. In
the second step, for each domino in a good pair we move one unit step,
upwards if it is N, downwards  if it is S, to the left if it is W and
to the right if it is E. After these two steps what remains to
completely fill $A_{k+1}$ are $2\times 2\,$-blocks. In the third and
last step we fill each $2\times 2\,$-block with a vertical pair with
probability $q$ and with a horizontal pair with probability
$1-q$. This generates a random tiling of the Aztec diamond $A_{k+1}$ with the
probability (\ref{2.1}), where $q=w^2(1+w^2)^{-1}$. 

The shuffling
algorithm translates into a PNG-type growth procedure for the cascade
of height functions $\{h_k(x,n)\}_{1\le k\le n}$ defined above.  We
will now define the discrete PNG-type growth model which we obtain. We
have shifted the picture 1/2 unit upwards compared to the one we
obtained from the modified DR-paths above. The ``height paths''
are built from {\it plus steps}, which go from $(m,n)$ to $(m+1,n+1)$ and
consists of straight line segments between the points $(m,n)$,
$(m+1/2,n)$, $(m+1/2,n+1)$ and $(m+1,n+1)$, {\it minus steps}, which
go from $(m,n)$ to $(m+1,n-1)$ and consists of 
straight line segments between the points $(m,n)$,
$(m+1/2,n)$, $(m+1/2,n-1)$ and $(m+1,n-1)$. Finally, we have {\it zero
  steps}, which are line segments from $(m,n)$ to $(m+2,n)$. The
initial configuration, $t=0$, has zero steps between $(n,0)$ and
$(n+2,0)$, $n\in 2\mathbb{Z}+1$. At time $t=k$, $k\ge 0$ we do the
following:

\begin{itemize}
\item[(i)] remove all zero steps;
\item[(ii)] move all plus steps one unit to the left and all minus
  steps one unit to the right;
\item[(iii)] if a plus step and a minus step pass each other in step
  (ii) they are removed;
\item[(iv)] add zero steps so that we obtain a connected curve from
  $(-\infty,0)$ to $(\infty,0)$;
\item[(v)] replace each zero step between $-(k+1)$ and $k+1$ with a
  combined plus and minus step independently with probability $q$.
\end{itemize}

We can define a {\it cascade of height curves} as follows. The level-$m$
curve initially has just zero steps between $(n,-(m-1))$ and
$(n+2,-(m-1))$, $n\in 2\mathbb{Z}+1$. At each time step we apply the
discrete PNG growth procedure independently to the levels $1,2,\dots$
with the condition that the level-$m$ curve cannot touch or intersect
the level-$(m+1)$ curve, $m\ge 1$. If that happens in the random
growth step, then this growth event is suppressed. Note har only a
finite number of levels are changed at time $k$. The shuffling
procedure induces an evolution of the modified DR-paths, and this is
exactly the cascade of PNG growth models just defined. A plus step
corresponds to a W-domino, a minus step to an E-domino and a zero step
to an S-domino. That the level-1 modified DR-paths evolves exactly
according to the PNG-growth rules is immediate by comparison. All
dominoes above the level-1 DR-path are N-dominoes and move upwards one
step. Removing all bad vertical pairs corresponds to the annihilation
(step (iii)) in the PNG-growth rule. Filling in $2\times 2\,$-blocks
is exactly the random growth, step (v) in the PNG growth rule. A
somewhat more elaborate argument shows that the whole shuffling
algorithm corresponds to the cascade defined above. We will not give
the details. 

Consider the level-1 DR-paths of type I, i.e. the upmost one. In CS-I
it goes from $(1,0)$ to $(n+1,n)$ through the points $(i_k,j_k)$,
$1\le k\le p$, where $(i_0,j_0)=(1,0)$, $(i_p.j_p)=(n+1,n)$ and
$(i_{k+1},j_{k+1})-(i_k,j_k)=(1,0)$, $(0,1)$ or $(1,1)$. The dominoes
in the north polar zone are the dominoes above this path. Set
\begin{equation}\label{2.60'}
n-\lambda_\ell=\max\{ j_k\,;\,i_k=\ell\},\quad 1\le\ell\le n+1.
\end{equation}
Then $\lambda=(\lambda_1,\dots,\lambda_{n+1})$ is a partition and this
is the partition associated with the north polar zone in
\cite{JPS}. 
(We can see this partition by marking each domino in the north polar
zone with a point at the center. These points will lie in the integer
lattice $\mathcal{L}_I$ in CS-I.)
Set
\begin{equation}
\Lambda (n)=\{(i,j)\in\mathbb{Z}_+^2\,;\,1\le j\le\lambda_i,\,1\le
i\le n+1\},\notag
\end{equation}
which is a random subset of $\mathbb{Z}_+^2$. 

Let
$w(i,j))_{(i,j)\in\mathbb{Z}_+^2}$ be independent geometric random
variables with parameter $q$, $P[w(i,j)=k]=(1-q)q^k$, $k\ge 0$, and
define
\begin{equation}\label{2.61}
G(M,N)=\max_{\pi}\sum_{(i,j)\in\pi} w(i,j),
\end{equation}
where the maximum is over all up/right paths from $(1,1)$ to $(M,N)$,
see \cite{Jo1}. Set
\begin{equation}
\Omega (n)=\{(i,j)\in\mathbb{Z}_+^2\,;\,G(M,N)+M+N-1\le n\}.
\end{equation}
We get $\Omega(n+1)$ from $\Omega(n)$ by independently
adding a point with probability $p=1-q$ to every corner in
$\Omega(n)$, see \cite{Jo1}, and because of this we call it the {\it corner
growth model}.
It is proved in \cite{JPS} that if the probability on
$\mathcal{T}(A_n)$ is defined by (\ref{2.1}) with $q=w^2(1+w^2)^{-1}$,
then the random sets $\Lambda(n)$ and $\Omega(n)$ have the same
distribution. 

Let $h(x,n)$ be the level-1 height function in the PNG growth model
defined above. We can relate the distribution function for this height
variable to the distribution function of $G(M,N)$. Let $P_0=\frac
12(n+2,n)$ be the midpoint in CS-I of the line segment from $(1,0)$ to
$(n+1,n)$ and set $P_k=P_0-k(1/2,1/2)$, $k=0,\pm 1,\dots,\pm
n$. Assume that both $k$ and $n$ are even. Then $h(k,n)\le 2m-1$ if
and only if $P_k+m(-1,1)$ is a point above the level-1 height curve,
which happens if and only if
\begin{equation}
(\frac{n-k}2-m+1,\frac {n+k}2-m+1)\in\Omega(n).
\notag
\end{equation}
Consequently, for $k$ and $n$ even,
\begin{equation}\label{2.62}
P[h(k,n)\le 2m-1]=P[G(\frac{n-k}2-m+1,\frac {n+k}2-m+1)\le 2m-1].
\end{equation}
Using this relation and the asymptotic results for $G(M,N)$ in
\cite{Jo1} we can show that the fluctuations of the height (and hence
of the temperate zone since the height describes the boundary of the
temperate zone) are of order $n^{1/3}$ and the appropriately rescaled
fluctuations converges to the Tracy-Widom distribution (\ref{1.11}).

As we have seen above the DR-paths can also be related to the zig-zag
particle configurations. Using this we can relate the distribution
function for $G(M,N)$ to the distribution of the rightmost particle in
the Krawtchouk ensemble.
\begin{lemma}\label{lem2.9}
If $K=t+N+M-1$, then
\begin{equation}\label{yy}
P[G(M,N)\le t]=P_{\text{Kr\,},M,K,q}[\max_{1\le j\le M} h_j\le t+M-1].
\notag
\end{equation}
\end{lemma}
\begin{Proof}
Considerthe zig-zag path $Z_r$ as defined above. It maps to the
zig-zag particle configuration $(h_1,\dots,h_r)$ with
$h_1<\dots<h_r$. The relation to the partition $\lambda$ defined by
(\ref{2.60'}) is that $\lambda_r=n-h_r$. Hence, $G(r,x)+r+x-1\le n$,
i.e. $\lambda_r\ge x$ in the tiling of $A_n$, if and only if $h_r\le
n-x$. Thus, by theorem \ref{thm2.2},
\begin{equation}
P[G(r,x)+r+x-1\le n]=P_{\text{Kr\,},M,K,q}[\max_{1\le j\le M} h_j\le
n-x],
\notag
\end{equation}
and this translates into (\ref{yy}).
\end{Proof}

We can now apply the edge scaling result for the Krawtchouk
ensemble. Using the integral formula for the Krawtchouk polynomials
(\ref{2.6}) and proceeding in the same way as in \cite{Jo1} for the
Meixner polynomials, we can prove that if $pt<q(1-t)$, $M=[Kt]$,
$0<t<1$, then
\begin{equation}\label{2.64}
\lim_{K\to\infty}P_{\text{Kr\,},M,K,q}[\max_{1\le j\le M} h_j\le
K\beta(t)+\xi\rho(t)K^{1/3}]=F(\xi)
\end{equation}
for each $\xi\in\mathbb{R}$. Here $F(\xi)$ is given by (\ref{1.11}) and
\begin{align}
\beta(t)&=(1-t)p+tq+2\sqrt{pqt(1-t)},\notag\\
\rho(t)&=\left(\frac{pq}{t(1-t)}\right)^{1/6}
(\sqrt{p(1-t)}+\sqrt{qt})^{2/3}
(\sqrt{q(1-t)}-\sqrt{pt})^{2/3}.\notag
\end{align}
We can now combine (\ref{yy}) and (\ref{2.64}) (allowing a somewhat
more complicated relation between $M,K$ and $t$) to give a new proof
of theorem 1.2 in \cite{Jo1}. Note that in the derivation of
(\ref{yy}) we have {\it not} used the RSK-correspondence which was central
to the approach in \cite{Jo1}. 

Let $L(\alpha)$ denote the length of a
longest increasing subsequence in a random permutation $\sigma$ from
$S_N$ where $N$ is a Poisson random variable with mean $\alpha$. Then
\begin{equation}
P[L(\alpha)\le n]=\lim_{N\to\infty}P[G(N,N)\le n]
\notag
\end{equation}
if we take $q=\alpha/N^2$, see \cite{Jo2}. Now, by (\ref{yy}), 
\begin{align}\label{2.65}
P[L(\alpha)\le n]&=\lim_{N\to\infty}
P_{\text{Kr\,},N,n+2N-1,\alpha/N^2}[\max_{1\le j\le M} h_j\le
n+N-1]\notag\\
&=\lim_{N\to\infty}\det(I-K_{\text{Kr\,},N,n+2N-1,\alpha/N^2})_{\ell^2
  (\{n+N,\dots,n+2N-1\})}.
\end{align}
Using the formulas (\ref{2.6}) and (\ref{2.7}), it follows that the
last expression of (\ref{2.65}) equals
\begin{equation}
\det(I-B_\alpha)_{\ell^2(\{n,n+1,\dots\})},
\notag
\end{equation}
where $B_\alpha$ is the {\it discrete Bessel kernel},
\begin{equation}
B_\alpha(x,y)=\sqrt{\alpha}
\frac{J_x(2\sqrt{\alpha})J_{y+1}(2\sqrt{\alpha})-
J_{x+1}(2\sqrt{\alpha})J_y(2\sqrt{\alpha})}{x-y},
\notag
\end{equation}
and we have rederived a result in \cite{BOO} and \cite{Jo2}. Precise
asymptotics for $L(\alpha)$ was first studied in \cite{BDJ}. We see
that the longest increasing subsequence problem can be found in a
limit of the Aztec diamond. In the same limit the discrete PNG model
defined above, appropriately rescaled, converges to the PNG model
studied in \cite{PrSp}.

\subsection{Zig-zag paths for domino tilings of the plane}

Consider the squares $(m,n)+[-1/2,1/2]^2$, where
$(m,n)\in\mathbb{Z}^2$. A domino tiling of the plane, which we
identify with $\mathbb{C}$, is a covering of the whole plane by
$2\times 1$ or $1\times 2$ rectangles, dominoes, where each domino
covers exactly two of the basic squares. We can also think of this as
a dimer configuration of the graph with vertices $(m,n)$ and edges
between nearest neighbour vertices. A domino covers the neighbouring
squares with centers $P$ and $Q$ if and only if the edge between $P$
and $Q$ are covered by a dimer. We will switch between the domino and
dimer languages whenever it is convenient. Colour the points $(m,n)$
with $m+n$ even black and the other points white, and give the
corresponding square the same colour. Consider the line $y=-x$. A
domino tiling of the plane induces an infinite zig-zag path around
black squares in complete analogy to the zig-zag paths in the Aztec
diamond. We can map the zig-zag path to a particle configuration by
saying that we have a particle at $x$ if and only if 
the zig-zag path goes from
$x-1/2-i(x-1/2)$ to $x+1/2-i(x+1/2)$ via $x+1/2+i(-x+1/2)$, i.e. an east-south
step. Note thst we have a particle at $x$ if and only if either the
edge from $x-1-ix$ to $x-ix$ or the edge from $x-i(x+1) $ to $x-ix$ is
covered by a dimer. In this way we get a particle configuration in
$\mathbb{Z}$.

There is a unique translation invariant measure $\mu$ of maximal
entropy, the {\it Burton-Pemantle measure}, on the space of domino tilings
of the plane, see \cite{BP} and \cite{Ke1}. This measure induces a
probability measure on zig-zag paths and hence on particle
configurations in $\mathbb{Z}$; we get a point process on
$\mathbb{Z}$. We want to show that this is a determinantal point
process, \cite{So2}, given by the discrete sine kernel. Let $E$ be a set of
disjoint edges, i.e. they do not share a vertex, in the
$\mathbb{Z}^2$-graph and let $U_E$ be the set of dimer configurations
which contain $E$. Let $P$ be a white vertex and give the edge between
$P$ and $P+z$ the {\it weight} $z$, where $z=\pm 1,\pm i$. Assume that
the edges in $E$ cover the black vertices $b_1,\dots,b_k$ and the
white vertices $w_1,\dots,w_k$. It is proved in \cite{Ke1},\cite{Ke2},
using techniques by Kasteleyn, \cite{K}, that
\begin{equation}
\mu(U_E)=a_E\det(P(b_i-w_j))_{i,j=1}^k,
\notag
\end{equation}
where $a_E$ is the product of the weights of the edges in $E$ and
\begin{equation}
P(x+iy)=\frac
1{4\pi^2}\int_{-\pi}^\pi\int_{-\pi}^\pi\frac{e^{i(x\theta-y\phi)}} 
{2i\sin\theta+2\sin\phi}d\theta d\phi.
\notag
\end{equation}
Using this we can prove the following result.
\begin{theorem}\label{thm2.9}
The probability of having particles at positions $x_1,\dots,x_m$ in
the zig-zag point process defined above is
\begin{equation}\label{2.66}
P[x]=\det\left(\frac{\sin\frac{\pi}2 (x_j-x_k)}{\pi
    (x_j-x_k)}\right)_{j,k=1}^m. 
\end{equation}
\end{theorem}
\begin{Proof}
If we have a particle at $x$, then one of the edges $x-ix,x-ix-1$ or
$x-ix,x-ix-i$ is covered by a dimer. We take $E_z$ to be the set of
edges $x_j-ix_j,x_j-ix_j-z_j$, where $z_j=1$ or
$=i$, $1\le j\le m$. 
Then, $a_E=z_1\dots z_m$, $b_j=x_j-ix_j$ and $w_k=x_k-ix_k-z_k$,
so that
\begin{equation}
P(b_j-w_k)=P(x_j-x_k-i(x_j-x_k)+z_k).
\notag
\end{equation}
Thus,
\begin{align}\label{2.67}
P[x]&=\sum_{z_j=1\,\,\text{or}\,\, i}
\mu(U_{E_z})=\sum_{z_j=1\,\,\text{or}\,\, i}z_1\dots z_m
\det( P(x_j-x_k-i(x_j-x_k)+z_k))\notag\\
&=\sum_{\sigma\in S_m}\text{sgn\,}(\sigma)\sum_{z_j=1\,\,\text{or}\,\, i}
\prod_{j=1}^mz_{\sigma(j)}P(x_j-x_{\sigma(j)}-i(x_j-
x_{\sigma(j)})+z_{\sigma(j)})\notag\\
&=\sum_{\sigma\in S_m}\text{sgn\,}
(\sigma)\prod_{j=1}^m K(x_j-x_{\sigma(j)})=
\det(K(x_j-x_k))_{j,k=1}^m,
\end{align}
where
\begin{align}
K(u)=\sum_{z_j=1\,\,\text{or}\,\, i}zP(u-iu+z)=
P(u+1-iu)+iP(u+(-u+1)i).\notag
\end{align}
A computation shows that $P(-y-ix)=i(-1)^yP(x+iy)$ and thus
\begin{equation}
P(u+1-iu)=i(-1)^{-u-1}P(u-i(u+1))=-i(-1)^uP(u-i(u+1)).\notag
\end{equation}
We obtain
\begin{align}
K(u)&=i(P(u+i(-u+1))-(-1)^uP(u-i(u+1))\notag\\
&=\frac
1{4\pi^2}\int_{-\pi}^\pi\int_{-\pi}^\pi\frac{e^{i(u\theta+(u-1)\phi)} 
-(-1)^ue^{i(u\theta+(u+1)\phi)}}{2i\sin\theta+\sin\phi}d\theta d\phi\notag\\
&=\frac 1{4\pi^2}\int_{-\pi}^\pi\int_{-\pi}^\pi\frac{e^{iu(\theta+\phi)}
(e^{-i\phi}-(-1)^ue^{i\phi})}{2i\sin\theta+\sin\phi}d\theta
d\phi\notag\\
&=\frac 1{4\pi^2}\int_{-\pi}^\pi\int_{-\pi}^\pi\frac{e^{i\theta
    u}(e^{-i\phi}-
(-1)^ue^{i\phi})}{2i\sin(\theta-\phi)+2\sin\phi}d\theta d\phi 
\notag
\end{align}
Set
\begin{equation}
G(\theta,u)=\frac
1{2\pi}\int_{-\pi}^\pi\frac{e^{-i\phi}-(-1)^ue^{i\phi}} 
{2i\sin(\theta-\phi)+2\sin\phi}d\phi
\notag
\end{equation}
so that
\begin{equation}
K(u)=\frac i{2\pi}\int_{-\pi}^\pi e^{i\theta u}G(\theta,u)d\theta.
\notag
\end{equation}
If we write
\begin{equation}
G(\theta,u)=\frac 1{2\pi i}\int_\gamma\frac{1-(-1)^uz^2}{e^{i\theta}
+i-(e^{-i\theta}+i)z^2}\frac{dz}z
\notag
\end{equation}
we can use residue calculus to see that
\begin{equation}
G(\theta,u)=
\begin{cases}
\frac{(-1)^u}{e^{-i\theta}+i} &\text{if $-\pi<\theta<0$}\\
\frac{1}{e^{i\theta}+i} &\text{if $0<\theta<\pi$}
\end{cases}.
\end{equation}
Thus,
\begin{align}
K(u)&=\frac
i{2\pi}\int_{-\pi}^0\frac{(-1)^ue^{iu\theta}}{e^{-i\theta}+i}d\theta 
+\frac i{2\pi}\int_0^\pi\frac{e^{iu\theta}}{e^{i\theta}+i} \notag\\
&=i\frac{1-(-1)^u}{2\pi u}=\frac{\sin\frac{\pi u}2}{\pi u}
(-1)^{u/2}.\notag
\end{align}
Inserting this in (\ref{2.67}) gives
\begin{equation}
p[x]=\det\left(\frac{\sin\frac{\pi}{2}(x_j-x_k)}{\pi
    (x_j-x_k)}(-1)^{\frac {x_j-x_k}2}\right)_{j,k=1}^m=
\det\left(\frac{\sin\frac{\pi}{2}(x_j-x_k)}{\pi
    (x_j-x_k)}\right)_{j,k=1}^m.
\notag
\end{equation}
\end{Proof}
If we consider the zig-zag path through the center of the Aztec
diamond, the measure on the zig-zag particle configurations has
determinantal correlation functions, (\ref{2.6'}), by
theorem \ref{thm2.2}. Take $r=n/2$ and $p=q=1/2$. In this case the equilibrium
measure, (\ref{2.22'}), has density $\rho'(x)=1/2$, $0\le x\le 1$. Hence, we
can use lemma \ref{lem2.8} to show that the limiting point pricess, as
$n\to\infty$, has determinantal correlation functions with kernel
\begin{equation}
K(x,y)=\frac{\sin\frac{\pi}{2}(x-y)}{x-y},\notag
\end{equation}
i.e. exactly the same as in theorem \ref{thm2.9}. This is consistent with the
conjecture, \cite{JPS}, \cite{CEP} 
that in the center of the Aztec diamond a random tiling
looks like a tiling of the plane under the Burton-Pemantle measure.

\section{The Schur measure and non-intersecting paths}
In section 2 we obtained the distribution function for the
last-passage random variable $G(M,N)$, (\ref{2.61}), using the
non-intersecting paths in the Aztec diamond. It is natural to inquire
whether the Meixner ensemble which is used to study $G(M,N)$ in
\cite{Jo1} can also be obtained in a natural way using
non-intersecting paths. The picture will again be a cascade of
PNG-type growth models but different from the one studied in
sect. 2.4.

We will define a certain random growth model, or rather a cascade of growth
models, such that the probability distribution of the heights above
the origin is given by the Schur measure introduced by Okounkov,
\cite{Ok1}. The cascade of growth models is actually equivalent with the
Robinson-Schensted-Knuth correspondence. Viennot, [Vi], gave a geometric
construction of the RSK-correspondence for permutations, often called
the ``shadow'' construction, see also \cite{Sa}, sect. 3.8. A
permutation in $S_n$ can be described by putting $n$ points randomly
in the unit square, Hammersley's picture, \cite{Ha}. If one applies
the first step in Viennot's construction and interprets the paths
(shadow lines) as space-time paths one gets the polynuclear growth
model (PNG) as introduced in \cite{PrSp} by Pr\"ahofer and
Spohn. Thus we can equivalently think in terms of a growth model. If
one takes the full Viennot construction one gets a cascade of growth
models as proved by Okounkov, \cite{Ok2}. Okounkov did not base his
presentation on Viennot's construction, but instead used the
formulation of the RSK-correspondence in \cite{BF}. There is a
generalization of Viennot's construction to the case of an integer
matrix (generalized permutation) called the ``matrix ball
construction'', see \cite{Fu}. By the same argument, this can be
translated into a growth model and will lead to the Schur
measure. This growth model is also given in \cite{Ok2}. We will
present a somewhat modified version of this growth model, which avoids
the limiting procedure in \cite{Ok2}. The interesting thing is that
this version leads directly to families of non-intersecting paths,
namely the standard ones which can be used to obtain the Schur
polynomials, \cite{Sa}, \cite{St}.

Let $W=(w(i,j)_{i,j=1}^n$ be an $n\times n$-matrix with non-negative
integer elements, and set $w(i,j)=0$ if $n\notin\{1,\dots,n\}^2$. The
integer-valued 
{\it height functions} in the cascade of growth models are denoted by
$h_k(x,t)$, $1\le k\le n$, where $h_k(x,t)$ is the height above
$x\in\mathbb{R}$ at time $t\in\mathbb{N}$ of the {\it level k
  growth process}. The height curves $x\to h_k(x,t)$ do not intersect,
$h_k(x,t)-h_{k+1}(x,t)\ge 1$, $1\le k\le n$ for all $x$ and $t$. The
height curves will grow by the addition of unit squares and the growth
procedure is defined as follows.

The vertical ``sides'' of $h_k(x,t)$ will be labelled. We can think of
the curve $x\to h_k(x,t)$ as a lattice path starting at $(-k+1,
-2n+1/2)$, ending at $(-k+1, 2n-1/2)$ and taking unit steps up, to the
right or down. Each unit step up is labelled by $a_j$ for some $j$,
$1\le j\le n$ and each unit step down is labelled by $b_k$ for some
$k$, $1\le k\le n$. Call these vertical sides {{\it left} and {\it 
right vertical sides} respectively. At time 0, $h_k(x,0)=
-(k-1)$, $1\le k\le n$ and there are no vertical sides. 

Assume that $h_k(x,t)$ has been defined for $t\le m-1$, $1\le k\le n$
with labels on the vertical sides and such that the distance between a
left vertical and a right vertical side is always odd. Furthermore 
$h_k(x,t)-h_{k+1}(x,t)\ge 1$. We will define $h_k(x,m), 1\le k\le n$
so that it has the same properties. For each $x\in\mathbb{Z}$,
$n\in\mathbb{N}$ and $1\le k\le n$, we have a set
$\mathcal{B}^{(k-1)}(x,t)$ of unit squares with labelled vertical
sides. $\mathcal{B}^{(0)}(x,t)$ contains $w(i,j)$ unit squares with
the left vertical side labelled $a_i$ and the right vertical side
labelled $b_j$, where $i=(t+x+1)/2$, $j=(t-x+1)/2$. Recall that
$w(i,j)=0$ if $n\notin\{1,\dots,n\}^2$ so that
$\mathcal{B}^{(0)}(x,t)$ is empty for odd $x$ at odd times $t$, and for even
$x$ at even times $t$. The $\mathcal{B}^{(k-1)}(x,t)$, $k\ge 2$ are defined
recursively in the growth procedure. Assume that
$\mathcal{B}^{(\ell -1)}(x,m)$ has been defined for some $\ell$, $1\le
\ell\le n$. We will define $h_\ell(x,m)$ and
$\mathcal{B}^{(\ell)}(x,m)$.

({\bf Horizontal growth}). Move each left vertical side in
$h_\ell(x,m-1)$ one unit to the left, and each right vertical side one
unit to the right, fill in with horizontal segments so that we get a
connected curve,  and denote the resulting height function by
$h_\ell^\ast(x, m-1)$. The labels move together with the vertical
sides. Since the distance between a left and a right vertical side is
odd they cannot meet at the same point. Set $u=h_\ell(x-1,m-1)-h_\ell(x,m-1)$
and $v=h_\ell(x+1,m-1)-h_\ell(x,m-1)$. If $z=\min(u,v)>0$, then a
right vertical side, with labels $b_{r_1},\dots,b_{r_u}$ (ordered in
the upwards direction), will cross a left vertical side, with labels
$a_{s_1},\dotsb_{s_v}$. If this happens we put $z$ unit squares in
$\mathcal{B}^{(\ell)}(x,m)$ with vertical sides labelled $a_{r_j},
b_{s_j}$, $1\le j\le z$. This defines $\mathcal{B}^{(\ell)}(x,m)$.

({\bf Vertical growth}). Next we put the labelled squares in
$\mathcal{B}^{(\ell -1)}(x,m)$ on top of $h_\ell^\ast(x, m-1)$ at
$x$ for all $x$, i.e. between $x-1/2$ and $x+1/2$. The result is
$h_\ell(x,m)$. Note that all vertical sides are labelled, and that by the
way that the $\mathcal{B}^{(\ell -1)}(x,m)$:s were defined, the
distance between left and right vertical sides is always odd.

In this way we recursively define $h_k(x,t)$, with labelled vertical
sides for $1\le k,t\le n$, $x\in\mathbb{R}$.

Next we want to describe the final configuration. 
Let $\lambda=(\lambda_1,\dots,\lambda_n)$ be a partition and let
$\mathcal{P}(\lambda;a)$ denote the set of all non-intersecting,
labelled, up/right lattice paths $\Gamma=\{\Gamma_k\}_{k=1}^n$, where
$\Gamma_k$ starts at $(-k+1,-2n+1/2)$ and ends at
$(-1/2,\lambda_j-j+1)$ and where all up-steps have an $x$-coordinate of
the form $2(j-n)-1/2$ for some $j$, $1\le j\le n$. Each unit length of
the vertical sides with the $x$-coordinate $2(j-n)-1/2$ is labelled
by $a_j$, $\le j\le n$. Let $\tilde{\mathcal{P}}(\lambda;b)$ be the
paths we obtain by reflecting the paths in $\mathcal{P}(\lambda;a)$ 
in the $y$-axis and putting the label
$b_k$ on each unit vertical side with $x$-coordinate $2(n-k)+1/2$,
$1\le k\le n$. If we join a family of non-intersecting paths 
from $\mathcal{P}(\lambda;a)$ to a family from
$\tilde{\mathcal{P}}(\lambda;b)$ by adding horizontal segments from
$(-1/2,\lambda_j-j+1)$ to $(1/2,\lambda_j-j+1)$, $1\le j\le n$, we
obtain a family of non-intersecting, labelled height curves. Let
$\mathcal{H}(\lambda;a,b)$ denote all the families of non-intersecting
height curves obtained in this way.

We claim that $H=\{h_k(x,2n-1)\}_{k=1}^n$ constructed as above belongs
to $\mathcal{H}(\lambda;a,b)$ if we put
$\lambda_j=h_j(0,2n-1)+j-1$. We write $H=(\Gamma,\tilde{\Gamma})$,
where $\Gamma\in\mathcal{P}(\lambda;a)$ and $\tilde{\Gamma}\in
\tilde{\mathcal{P}}(\lambda;b)$. To see this note that a square with
labels $a_j,b_k$ is introduced at level 1 at position $j-k$ at time
$j+k-1$. In the remaining time $2n-1-(j+k-1)=2n-j-k$, the left
vertical side moves $2n-j-k$ steps to the left and the right vertical
side $2n-j-k$ steps to the right. Thus the $a_j$-label ends up at a
position with $x$-coordinate $j-k-1/2-(2n-j-k)=2(j-n)-1/2$, and the
$b_k$-label ends up at  a position with $x$-coordinate
$j-k+1/2+2n-j-k=2(n-k)+1/2$. Thus all left (right) vertical sides end
up to the left(right) of the origin at the correct positions. Thus, we
obtain a map from the set of $n\times n$ integer matrices to 
$\cup_{\lambda}\mathcal{H}(\lambda;a,b)$, where the union is over all
partitions $\lambda=(\lambda_1,\dots,\lambda_n)$. This map is
invertible, and hence we obtain a bijection. That the map is
invertible follows from the fact that the growth procedure can be
reversed! We start at the bottom level, $h_n(x,t)$, but now right
vertical sides move to the left and left vertical sides to the
right. If they cross we move squares up to the next level, just as we
moved them down before. At the upper level they tell us where we
should split and introduce new left and right vertical sides. Those
squares that are taken out at the top level give us the entries in the
matrix. If we take out $m$ squares with labels $a_j,b_k$ we put the
number $m$ at position $(j,k)$ in the matrix. (We can also find the
position from the time $t$ when and the place $x$ where the squares
are removed, $x=j-k$, $t=j+k-1$.) We will not describe all the
details of this reverse procedure. In a sense the cascade of growth
models records the history of the growth process. When two vertical
segments pass each other at a certain level information is lost at
this level, but this information is recorded at the next level.

If the $w(j,k)$:s are independent geometric random variables with
$P[w(j,k)=m]=(1-a_jb_k)(a_jb_k)^m$, $m\ge 0$, then the probability of
a particular integer matix $W=(w(j,k))_{j,l=1}^n$ is
\begin{equation}
\prod_{j,k=1}^n(1-a_jb_k)\prod_{j,k=1}^n (a_jb_k)^{w(j,k)}=
\prod_{j,k=1}^n(1-a_jb_k)\omega(W),\notag
\end{equation}
where
\begin{equation}
\omega(W)=\prod_{j=1}^na_j^{\sum_k w(j,k)}\prod_{k=1}^nb_k^{\sum_j
  w(j,k)}.\notag
\end{equation}
If we interpret the labels on the vertical sides of
$\mathcal{P}(\lambda;a)$ as weights and define the weight,
$\omega(\Gamma)$, of an element $\Gamma\in\mathcal{P}(\lambda;a)$ as
the product of the weights of all vertical sides (horizontal sides
have weight 1), then we see that the growth procedure defined above
transports the weights in the correct way; if $W$ has weight
$\omega(W)$, and $W$ maps to $(\Gamma,\tilde{\Gamma})$ in
$\mathcal{H}(\lambda;a,b)$, then
\begin{equation}
\omega(W)=\omega(\Gamma)\omega(\tilde{\Gamma}).\notag
\end{equation}
The total weight of all up/right paths from $(-k+1,-2n+1/2)$ to
$(-1/2,\lambda_j-j+1)$ where all vertical sides have $x$-coordinates
of the form $2(i-n)-1/2$, $1\le i\le n$, is given by
\begin{equation}
\sum_{m_1+\dots+m_n=\lambda_j-j+k}a_1^{m_1}\dots a_n^{m_n}=
h_{\lambda_j-j+k}(a_1,\dots,a_n);\notag
\end{equation}
We have $m_i$ vertical steps with $x$-coordinate $2(i-n)-1/2$ and
these have weight $a_i$. Here $h_m(a_1,\dots,a_n)$ is the complete
symmetric polynomial of degree $m$ in $n$ variables; $h_m(a)\equiv 0$
if $m<0$.

The Lindstr\"om-Gessel-Viennot method, which is the discrete analogue
of the Karlin-McGregor theorem, see for example \cite{Stem}, p. 98 for a
precise statement, gives
\begin{equation}\label{LGV}
\sum_{\Gamma\in\mathcal{P}(\lambda;a)}\omega(\Gamma)=
\det(h_{\lambda_j-j+k}(a))_{j,k=1}^N=s_\lambda(a).
\end{equation}
The expresion in the middle can be taken as the definition of the
{\it Schur polynomial} $s_\lambda(a)$. The growth procedure above defines a
map $S(W)=\lambda$ from the integer matrix $W$ to the partition
$\lambda$ defined by the succesive heights. We obtain
\begin{align}\label{sch}
P[S(W)=\lambda]&=\sum_{W:S(W)=\lambda}\prod_{j,k=1}^n(1-a_jb_k)\omega(W)\\
&=\prod_{j,k=1}^n(1-a_jb_k)\sum_{\Gamma\in\mathcal{P}(\lambda;a),
\tilde{\Gamma}\in\tilde{\mathcal{P}}(\lambda;b)}\omega(\Gamma)
\omega(\tilde{\Gamma})\notag \\ 
&=\left[\prod_{j,k=1}^n(1-a_jb_k)\right]s_\lambda(a)s_\lambda(b)\doteq
P_{\text{Schur}}[\lambda],\notag
\end{align}
the {\it Schur measure} on partitions, introduced in \cite{Ok1}.

\begin{remark}
\rm Note that an element $\Gamma=\{\Gamma_k\}_{k=1}^n$ in
$\mathcal{P}(\lambda;a)$ corresponds to a unique semi-standard
tableaux $T$ with shape $\lambda$, $\text{sh\,}(T)=\lambda$. If
$\Gamma_k$ has $r_j$ vertical steps at  $2(j-n)-1/2$, we put $r_j$
$j$:s, $j=1,\dots, n$ in weakly increasing order in the $\lambda_k$
boxes in row $k$. Similarly an element in
$\tilde{\mathcal{P}}(\lambda;b)$ gives a semistandard tableaux
$D$. Thus we obtain a one-to-one map $W\to (S,T)$, which is exactly
the RSK-correspondence. If we let $m_j(T)$ denote the number of $j$:s
in $T$, we see that we also have, by (\ref{LGV}), 
\begin{equation}
s_\lambda(a)=\sum_{\Gamma\in\mathcal{P}(\lambda;a)}
\omega(\Gamma)=\sum_{T:\text{sh\,}=\lambda}a_1^{m_1(T)}\dots
a_n^{m_n(T)},
\notag
\end{equation}
which is the combinatorial definition of the Schur polynomial, see
\cite{St}.\it
\end{remark} 

Consider now the random variable $G(M,N)$ defined by (\ref{2.61}).
From its definition it is clear that we can compute $G(M,N)$
recursively by
\begin{equation}\label{3.0}
G(M,N)=\max(G(M-1,N),G(M,N-1))+w(M,N).
\end{equation}
We want to show that, for $1\le M,N\le n$,
\begin{equation}\label{3.1}
G(M,N)=h_1(M-N,M+N-1),
\end{equation}
in particular $G(n,n)=h_1(0, 2n-1)=\lambda_1$, which is a well known
property of the RSK-correspondence, \cite{Sa}, \cite{Kn}.
It follows from the definition of the growth process above that
\begin{equation}\label{3.2}
h(x,t)=\max(h_1(x-1,t-1),h(x,t-1), h(x+1, t-1))+w(\frac{t+x+1}2,
\frac{t-x+1}2).
\end{equation}
(Recall that $w(j,k)=0$ if $(j,k)\notin\{1,\dots, n\}^2$.) 
Write $\tilde{G}(j,k)=h(j-k,j+k-1)$. Then (\ref{3.2}) becomes,
\begin{equation}\label{3.3}
\tilde{G}(j,k)=\max(\tilde{G}(j-1,k),\tilde{G}(j-\frac 12,k-\frac 12), 
\tilde{G}(j,k-1))+w(j,k).
\end{equation}
We will now use induction. If $M+N-1=1$, then $M=N=1$ and 
$G(1,1)=w(1,1)=\tilde{G}(1,1)$. Assume that $G(M,N)=\tilde{G}(M,N)$
for $M+N-1<k$. We want to show that it is true for $M+N-1=k$. Note
that, by (\ref{3.3}),
\begin{equation}\label{3.4}
 \tilde{G}(M-\frac 12,N-\frac 12)=\max(\tilde{G}(M-\frac 32,N-\frac
 12), \tilde{G}(M-1,N-1), \tilde{G}(M-\frac 12, N-\frac 32))
\end{equation}
since $w(M-1/2,N-1/2)=0$. By our assumption and (\ref{3.3}),
\begin{equation}\label{3.5}
G(M-1, N)=\tilde{G}(M-1,N)\ge \tilde{G}(M-\frac 32, N-\frac 12)
\end{equation}
and
\begin{equation}\label{3.6}
G(M, N-1)=\tilde{G}(M,N-1)\ge \tilde{G}(M-\frac 12, N-\frac 32).
\end{equation}
Furthermore, by our assumption, $\tilde{G}(M-1,N-1)=G(M-1, N-1)$ and
$G(M-1,N)\ge G(M-1, N-1)$ by (\ref{3.0}). Combining this with
(\ref{3.4}), (\ref{3.5}) and (\ref{3.6}) we see that
\begin{equation}\label{3.7}
\max(G(M-1,N),G(M, N-1))\ge \tilde{G}(M-\frac 12, N-\frac 12).
\end{equation}
Consequently, by (\ref{3.3}), our assumption and (\ref{3.7})
\begin{align}
\tilde{G}(M,N)&=\max(G(M-1,N),\tilde{G}(M-\frac 12,N-\frac 12),
G(M,N-1)) +w(M,N)\\ \notag
&=\max(G(M-1,N),G(M,N-1))+w(M,N)=G(M,N),\notag
\end{align}
which completes the proof.

From (\ref{3.1}) and (\ref{sch}) we obtain
\begin{equation}
P[G(M,N)\le t]=[\prod_{j,k=1}^n(1-a_jb_k)]
\sum_{\lambda:\lambda_1\le t}s_\lambda(a)s_\lambda(b).
\notag
\end{equation}
Using the fact that the Schur measure has determinantal correlation
functions, which was proved in \cite{Ok1}, see also 
\cite{Jo3}, we see that this equals a
Fredholm determinant with a certain kernel, and this can be exploited
for the asymptotic analysis.

\begin{remark}
\rm The case when $w(i,j)$ are independent exponential random
variables with mean 1 can be obtained as a limit of the geometric case
as discussed in \cite{Jo1}. We can take the same limit in the
construction above and this leads to a continuous analogue of the
RSK-correspondence. The resulting picture of paths can be viewed as
two families of $n$ non-intersecting Poisson processes with rate
$1/2$. (We take $a_j=1-1/2L$, $b_j=1-1/2L$, so that $a_jb_k\approx
1-1/L$ when $L$ is large. This gives the rate $1/2$ for the limiting
Poisson processes on both sides.) Take the vertical axis
in the negative direction as time axis, and the horizontal axis as
counting the number of events in the Poisson process. If the heights
above the origin in the cascade are $t_1>\dots>t_n$, then the $k$:th
process $A_k$ on the left (and right) start at 1 at the time $-t_k$
and end at $n+1-k$ at time 0. The probability for $X_i$ to go from 0
to $n-j$ in a time interval of length $t_i$ is
\begin{equation}\label{3.12}
P[X_i(t_i)=n-j]=e^{-t_i/2}\frac{(t_i/2)^{n-j}}{(n-j)!}.
\end{equation}
The Karlin/McGregor theorem can be generalized to unequal starting
times, see \cite{Karl2}, and we find that the probability for the $n$
non-intersecting paths to the left with the specified initial and
final positions is
\begin{equation}\label{3.13}
\det(e^{-t_i/2}\frac{(t_i/2)^{n-j}}{(n-j)!})_{i,j=1}^n=2^{-n(n+1)/2}
\prod_{j=1}^n\frac 1{j!}\Delta_n(t)\prod_{j=1}^ne^{-t_i/2}.
\end{equation}
The possible heights lie in $[0,\infty)$, so we obtain the probability
density
\begin{equation}\label{3.14}
\frac{1}{Z_n}\Delta_n(t)^2\prod_{j=1}^ne^{-t_i},
\end{equation}
where
\begin{equation}
Z_n=\int_{[0,\infty)^n}\Delta_n(t)^2\prod_{j=1}^ne^{-t_i}d^nt.
\notag
\end{equation}
Hence we have rederived the result of proposition 1.4 in \cite{Jo1} (in
the case $M=N$) using non-intersecting paths. The probability density
(\ref{3.14}) is a special case of the Laguerre ensemble. It is also
possible to consider the case when $w(i,j)$ is exponential with
parameter $\alpha_i+\beta_j$. This leads to a continuous analogue of
the Schur measure, which can also be obtained as a limit of the Schur
measure defined above with $a_i=1-\alpha_i/L$, $b_i=1-\beta_i/L$ as
$L\to\infty$. Using the methods of \cite{Jo3} we can derive the
correlation functions for the continuous Schur measure, and we can also
obtain proposition 1.4 in \cite{Jo1} in the case $M\neq N$.
\it
\end{remark}
\section{Random walks and rhombus tilings of a hexagon}
\subsection{Derivation of the Hahn ensemble}

In this section we will explore the relation between non-intersecting
random walk paths and another tiling problem. We will consider random
tilings of an abc-hexagon with rhombi, see \cite{CLP}, which are
directly related to so called boxed plane partitions, \cite{St}. An 
{\it abc-hexagon} has sides $a,b,c,a,b,c$ (in clockwise order) and
equal angles. We want to tile this region with unit rhombi (often
called lozenges) with angles $\pi /3$ and $2\pi/3$. The number of
possible tilings is given by MacMahon's formula, \cite{St},
\begin{equation}\label{4.1}
N(a,b,c)=\prod_{i=1}^a\prod_{i=1}^b\prod_{i=1}^c\frac{i+j+k-1}
{i+j+k-2}.
\end{equation}
We obtain a random tiling by picking each tiling with equal
probability. It is well known that a 
rhombus tiling can be described by
non-intersecting random walk paths, see for example \cite{Fisch}
and below.
This is the approach we will adopt here. We will give a new derivation
of the Hahn ensemble introduced in \cite{Jo2}.

Let $a,b,c\ge 1$ be given integers. Take $\mathbf{e}=(0,1/2)$ and
$\mathbf{f}=(\sqrt{3}/2,0)$ as basis vectors in our coordinate
system; all coordinates will refer to this choice of basis
vectors. Our hexagon has corners at $P_1=(0,0)$, $P_2=(b,-b)$,
$P_3=(a+b,a-b)$, $P_4=(a+b,a-b+2c)$, $P_5=(a,a+2c)$ and
$P_6=(0,2c)$. Consider random walks $S^k(m)$, $0\le m\le a+b$,
starting at $(0,2(k-1))$ and ending at $(a+b,a-b+2(k-1))$, $1\le k\le
c$,
\begin{equation}
S^k(m)=(0,2k)+\sum_{j=1}^m(1,X_j^k),
\notag
\end{equation}
where $X_j^k=\pm 1$ are independent Bernoulli random variables taking each
value with probability $1/2$. Assume that these random walks are
non-intersecting. We restrict to the case $a\ge b$; the case $a\le b$
is completely analogous. Set
\begin{align}
\alpha_m&=
\begin{cases}
-m & \text{if $0\le m\le b$}\\
m-2b & \text{if $b\le m\le a$}\\
m-2b & \text{if $a\le m\le a+b$},
\end{cases}\notag\\
\beta_m&=
\begin{cases}
m+2(c-1) & \text{if $0\le m\le b$}\\
m+2(c-1) & \text{if $b\le m\le a$}\\
2a-m+2(c-1) & \text{if $a\le m\le a+b$},
\end{cases}\notag
\end{align}
and $\gamma_m=(\beta_m-\alpha_m)/2$. Then $\alpha_m\le S^k(m)\le
\beta_m$.
Set $x_k=(S^k(m)-\alpha_m)/2$ and note that $0\le x_k\le \gamma_m$,
$1\le k\le c$. These numbers describe the points where the random
walks intersect a vertical line. Note the analogy with the previous
problems. Think of $x_1,\dots,x_c$ as the positions of $c$ {\it
  particles} in a discrete gas confined to $\{0,\dots,\gamma_m\}$. Let
$\xi_1<\dots<\xi_{L_m}$,
\begin{equation}\label{4.2}
L_m=\gamma_m+1-c=
\begin{cases}
m & \text{if $0\le m\le b$}\\
b & \text{if $b\le m\le a$}\\
a+b-m & \text{if $a\le m\le a+b$},
\end{cases}
\end{equation}
be the positions of the {\it holes}. The holes correspond to the
positions $(m,\alpha_m+2\xi_k)$ in our coordinate system.

There are three types of rhombi. {\it Type I} which are spanned by
$\mathbf{e}+\mathbf{f}$ and $2\mathbf{f}$, 
{\it Type II} which are spanned by
$\mathbf{e}-\mathbf{f}$ and $2\mathbf{f}$ and
{\it Type III} (called {\it vertical}) spanned by
$-\mathbf{e}+\mathbf{f}$ and $\mathbf{e}+\mathbf{f}$. (We will
sometimes call type I and II {\it horizontal}.) Given the
non-intersecting random walk paths we can now tile the hexagon with
rhombi as follows. If $X_{m+1}^k=1$ we put a type I rhombus at
$(m,S^k(m))$, and if $X_{m+1}=-1$ we put a type II rhombus at
$(m,S^k(m))$. Finally we put a type III rhombus at each of the hole
positions $(m,\alpha_m+2\xi_k)$. Note that the vertical rhombi are
associated with holes in the gas, whereas each particle is associated
with a horizontal rhombus. Our non-intersecting random walks
correspond to picking all of the possible hexagon tilings of the
abc-hexagon with equal probability. We want to compute the probability
measure induced on the particle/hole configurations on the vertical
line $x=m$. In order to be able to formulate the results 
we first define the Hahn and associated Hahn
ensembles. 

The {\it Hahn ensemble}, \cite{Jo2} is a probability measure on
$\{0,\dots, N\}^m$ defined by
\begin{equation}\label{4.3}
P_{N,m}^{(\alpha,\beta)}[h]=\frac
1{Z_{N,m}^{(\alpha,\beta)}}\Delta_m(h)^2
\prod_{j=1}^mw_N^{(\alpha,\beta)}(h_j),
\end{equation}
where $\alpha,\beta>-1$ are given parameters,
\begin{equation}\label{4.4}
w_N^{(\alpha,\beta)}(t)=\frac{(N+\alpha-t)!(\beta+t)!}{t!(N-t)!}
\end{equation}
a weight function and
\begin{equation}\label{4.5}
Z_{N,m}^{(\alpha,\beta)}=\sum_{h\in\{0,\dots, N\}^m}
\Delta_m(h)^2
\prod_{j=1}^mw_N^{(\alpha,\beta)}(h_j)
\end{equation}
a normalization constant. This ensemble is related to the Hahn
polynomials, \cite{NSU}, which are orthogonal on $\{0,\dots,N\}$ with
respect to the weight (\ref{4.4}). For some facts about these
polynomials see the proof of lemma \ref{lem4.2} below. Using the leading
coefficients of the normalized Hahn polynomials a standard
computation, \cite{Me}, gives
\begin{equation}\label{4.6}
Z_{N,m}^{(\alpha,\beta)}=m!\prod_{j=0}^{m-1}
\frac{j!(\alpha+j)!(\beta+j)!(\alpha+\beta+j+N+1)!(\alpha+\beta+j)!}
{(\alpha+\beta+2j)!(\alpha+\beta+2j+1)!(N-j)!}.
\end{equation}
The {\it associated Hahn ensemble} on $\{0,\dots, N\}^m$ is defined by
\begin{equation}\label{4.7}
\tilde{P}_{N,m}^{(\alpha,\beta)}[h]=\frac
1{\tilde{Z}_{N,m}^{(\alpha,\beta)}}\Delta_m(h)^2
\prod_{j=1}^m\tilde{w}_N^{(\alpha,\beta)}(h_j),
\end{equation}
where $\tilde{Z}_{N,m}^{(\alpha,\beta)}$ is a normalization constant
and 
\begin{equation}\label{4.8}
\tilde{w}_N^{(\alpha,\beta)}(t)=
\frac 1{t!(N-t)!(N+\alpha-t)!(\beta+t)!}.
\end{equation}
The ensembles (\ref{4.3}) and (\ref{4.7}) are related by a
particle/hole transformation, compare (\ref{2.13'}) above, and see the
proof of theorem \ref{thm4.1} below.

Let $\tilde{P}_m(x_1,\dots,x_c)$ denote the probability of having the
particle configuration $x_1,\dots,x_m$ at time $m$ (along the vertical
axis $x=m$), and $P(\xi_1,\dots,\xi_{L_m})$, $L_m$ given by
(\ref{4.2}), the probability of having the hole configuration
$\xi_1,\dots,\xi_{L_m}$ at time $m$.
\begin{theorem}\label{thm4.1}
If $a,b,c\ge 1$ are given integers, $a\ge b$, and we define
$a_m=|a-m|$, $b_m=|b-m|$, then
\begin{equation}\label{4.9}
\tilde{P}_m(x_1,\dots,x_c)=\tilde{P}_{\gamma_m,c}^{(a_m,b_m)}[x]
\end{equation}
and
\begin{equation}\label{4.10}
P_m(\xi_1,\dots,\xi_{L_m})=P_{\gamma_m,L_m}^{(a_m,b_m)}[\xi]
\end{equation}
\end{theorem}
\begin{Proof}
The number of random walk paths from $j$ to $k$ in $m$ steps is
($m+k-j$ even),
\begin{equation}
\binom{m}{\frac{m+k-j}2}=e_{\frac{m+k-j}2}(1^m),
\notag
\end{equation}
where $1^m=(1,\dots,1)\in\mathbb{N}^m$ and $e_n(x)$ is the elementary
symmetric function. We can now use the Karlin-McGregor,
Lindstr\"om-Gessel-Viennot argument to see that the number of
non-intersecting random walk paths from $(0,2(k-1))$ to
$(m,2x_k+\alpha_m)$ is
\begin{equation}\label{4.11}
A_m(x)\doteq\det(e_{\delta_m+x_k-j+1}(1^m))_{j,k=1}^c,
\end{equation}
where $\delta_m=(m+\alpha_m)/2$. Introduce the shifted particle
positions
\begin{equation}
s_k=x_{c+1-k}+\delta_m\in\{0,\dots,\gamma_m+\delta_m\},
\notag
\end{equation}
which means that we have introduced $\delta_m$ extra holes at the
positions $j-1$, $1\le j\le\delta_m$. Define the partition $\lambda$
by
\begin{equation}
\lambda_k=s_k+k-c,\quad 1\le k\le c.
\notag
\end{equation}
By reversing the order of the rows and columns in (\ref{4.11}) we
obtain
\begin{align}
A_m(x)&=\det(e_{\delta_m+x_{c+1-k}-c+j}(1^m))_{j,k=1}^c\notag\\
&=\det(e_{\lambda_k-k+j}(1^m))_{j,k=1}^c=s_{\lambda'}(1^m),
\end{align}
where $\lambda'$ is the conjugate partition to $\lambda$ and
$s_\lambda$ is the Schur polynomial, see (\ref{3.1}) above and \cite{Sa}.

If we set $r_j=c+j-1-\lambda_j'$, $1\le j\le L_m+\delta_m$, then
$\{r_1,\dots,r_{L_m+\delta_m}\}
\cup\{s_1,\dots,s_c\}=\{0,\dots,\gamma_m+\delta_m\}$ so the $r_k$ give the
positions of the (shifted) holes including the extra holes. The
positions of the original holes are given by
\begin{equation}\label{4.13}
\xi_j=r_{\delta_m+j}-\delta_m,\quad 1\le j\le L_m.
\end{equation}
Note that $L_m+\delta_m=m$ if $0\le m\le a$ and $L_m+\delta_m=a$ if
$a\le m\le a+b$. Let $\mu=\lambda'=(\lambda_1',\dots\lambda_m')$ if
$0\le m\le a$, and $\mu=(\lambda_1',\dots\lambda_m',0,\dots,0)$ ($m-a$
extra zeros) if $a<m\le a+b$. Then,
\begin{equation}\label{4.14}
A_m(x)=s_{\lambda'}(1^m)=s_\mu(1^m)=\prod_{1\le i\le j\le m}
\left(\frac{\mu_i-\mu_j+j-i}{j-i}\right), 
\end{equation}
by the classical formula for a Schur polynomial. We now want to
rewrite the right hand side of (\ref{4.14}) in terms of $\xi$ using
(\ref{2.13'}) and (\ref{4.13}). Some computation gives
\begin{equation}\label{4.15}
A_m(x)=C(m,a,b,c)\Delta_{L_m}(\xi),
\end{equation}
if $0\le m\le b$,
\begin{equation}\label{4.16}
A_m(x)=C(m,a,b,c)\prod_{j=1}^{L_m}\frac {(\xi_j+m-b)!}{\xi_j!} 
\Delta_{L_m}(\xi),
\end{equation}
if $b\le m\le a$, and
\begin{equation}\label{4.17}
A_m(x)=C(m,a,b,c)\prod_{j=1}^{L_m}\frac {(\xi_j+m-b)!(b+c-1-\xi_j)}
{\xi_j!(a+b+c-m-1-\xi_j)!} 
\Delta_{L_m}(\xi),
\end{equation}
if $a\le m\le a+b$, where $C(m,a,b,c)$ is a constant, e.g.
\begin{equation}\label{4.18}
C(m,a,b,c)=\prod_{j=0}^{L_m-1}\frac 1{(j+m-b)!}
\end{equation}
for $b\le m\le a$.

We will also write $W_m(\xi)=A_m(x)=$ the number of non-intersecting
random walks ending with hole configuration $\xi$ at time $m$. The
number of possible non-intersecting random walks coming from the right
side and going in the other direction is $W_{m'}(\xi')$, where
$\xi_j'=c+L_m-1-\xi_{L_m+1-j}$
and $m'=a+b-m$. This follows from the symmetry of the hexagon. The
total number of tilings, given $\xi$, at time $m$ is then
$W_m(\xi)W_{m'}(\xi')$. This can be computed by using (\ref{4.15}) -
(\ref{4.17}) with the result
\begin{equation}\label{4.19}
W_m(\xi)W_{m'}(\xi')=C^\ast(m,a,b,c)\Delta_{L_m}(\xi)^2\prod_{j=1}^{L_m}
w_{\gamma_m}^{(a_m,b_m)}(\xi_j),
\end{equation}
where $C^\ast(m,a,b,c)$ is a constant. Using (\ref{4.18}) we see that,
for $b\le m\le a$
\begin{equation}\label{4.20}
C^\ast(m,a,b,c)=C(m,a,b,c)C(m',a,b,c)=\prod_{j=0}^{b-1}\frac
1{(j+m-b)!(j+a-m)!}.
\end{equation}
It follows that the total number of tilings is
\begin{align}
N(a,b,c)&=C^\ast(m,a,b,c)\sum_{\xi\in\{0,\dots,\gamma_m\}^{L_m}}
\Delta_{L_m}(\xi)^2\prod_{j=1}^{L_m}
w_{\gamma_m}^{(a_m,b_m)}(\xi_j)\notag\\
&=C^\ast(m,a,b,c)Z_{\gamma_m,L_m}^{(a_m,b_m)},\notag
\end{align}
by (\ref{4.5}). This proves (\ref{4.10}). Note that by combining
(\ref{4.20}) and (\ref{4.6}) we obtain, after some computation (where
we take $m=b$),
\begin{equation}
N(a,b,c)=\prod_{j=0}^{b-1}\frac{j!(a+c+j)!}{(a+j)!(c+j)!},
\notag
\end{equation}
which proves MacMahon's formula (\ref{4.1}).

We now want to go from the variables $\xi$ to $x$. Since
$\{x_1,\dots,x_c\}\cup \{\xi_1,\dots,\xi_{L_m}\}=\{0,\dots,\gamma_m\}$
we can use (\ref{2.13'}) to get
\begin{equation}
\Delta_{L_m}(\xi)=\left(\prod_{j=1}^{\gamma_m}j!\right)
\left(\prod_{k=1}^c\frac 1{x_k!(\gamma_m-x_k)!}\right)\Delta_c(x).
\notag
\end{equation}
Using this it is straightforward to show that 
\begin{equation}
W_m(\xi)W_{m'}(\xi')=C_\ast(m,a,b,c) 
\Delta_c(x)^2\prod_{j=1}^c\tilde{w}_{\gamma_m}^{(a_m,b_m)}(x_j),
\notag
\end{equation}
and (\ref{4.9}) follows,
\end{Proof}

\subsection{Some asymptotic results}

Random tilings of a hexagon shows the same type of arctic ellipse
effect as the Aztec diamond, see .... We will have polar zones
associated with each corner of the hexagon. Consider the vertex
$P_6=(0,2c)$ of the hexagon. We will say that two type I rhombi are
{\it adjacent} if they share an edge. The {\it polar zone} associated
with $P_6$ is now defined as follows. If there is no type I rhombus
$R_0$ having $P_6$ as a vertex, the polar zone is empty. Otherwise a
type I rhombus $R$ belongs to the polar zone if there is a sequence of
rhombi $R_0,\dots,R_k=R$ such that $R_j$ and $R_{j+1}$ are
adjacent. Consider the horizontal rhombi immediately to the left of
the line $x=m$, the {\it m:th column}, there are $c$ of them. If
$Z_m=\xi_{L_m}=\max \xi_j$ is the position of the last hole,
i.e. vertical rhombus on the line $x=m$, then all the horizontal
rhombi above it in the $m\,$:th column are of type I and belong to the
polar zone of $P_6$. Hence, the boundary of this polar zone is
obtained by joining the points $A_1,B_1,\dots, A_a,B_a$, where
$A_m=(m-1,\alpha_m+2Z_m+1)$, $B_m=(m,\alpha_m+2Z_m)$ with straight
line segments. The boundary of this polar zone is thus related to the
position of the rightmost particle in the Hahn ensemble.

The asymptotic position of the rightmost particle and its large
deviation properties can be investigated using the results of
\cite{Jo1}, sect. 2.2. Consider the Hahn ensemble (\ref{4.3}). Assume
that $m/N\to t\in (0,1)$ and $\frac 1{N}(\alpha,\beta)\to
(\alpha_0,\beta_0)$ as $N\to\infty$. We have the limit
\begin{equation}
V(s)=-\lim_{N\to\infty}\frac 1m W_N^{(\alpha,\beta)}(ms)=-\frac 1t
U(ts), 
\notag
\end{equation}
where
\begin{align}
U(s)&=(1+\alpha_0-s)\log(1+\alpha_0-s)+(\beta_0+s)\log(\beta_0+s)\notag
\\&-s\log s   -(1-s)\log (1-s)-C.
\notag
\end{align}
Here we have introduced the modified weight,
\begin{equation}\label{4.20'}
W_N^{(\alpha,\beta)}(x)=\frac{\binom{\beta+x}{x}\binom{\alpha+N-x}{x}}
{\binom{\alpha+\beta+N+1}{N}},
\end{equation}
which only differs by a multiplicative constant. The {\it equilibrium
  measure}, $u_{\text{eq}}^{(t,\alpha_0,\beta_0)}(s)ds$, for the Hahn
  ensemble is the unique solution of the constrained, weighted
  variational problem

\begin{equation}\label{4.21}
F_V=\inf_{u\in\mathcal{A}}\left(\int_0^{1/t}\int_0^{1/t}
  \log|\sigma-s|^{-1} u(\sigma)u(s)d\sigma ds+\int_0^{1/t}
  V(s)u(s)ds\right), 
\end{equation}
where $\mathcal{A}=\{u\in
L^1[0,1/t]\,;\,\int_0^{1/t}u=1\,\,\text{and}\,\, 0\le u\le 1\}$. From
\cite{Jo1}, theorem 2.2, we obtain the following large deviation
result. Let $R=R(t,\alpha_0,\beta_0)$ be the right endpoint of the
support of the equilibrium measure, and let $\epsilon>0$ be
given. There are functions $L(R-\epsilon)$ and $J(R+\epsilon)$ such
that 
\begin{equation}\label{4.22}
\lim_{N\to\infty}\frac 1{N^2}\log P_{N,m}^{(\alpha,\beta)}[\frac
1m\max h_j\le R-\epsilon]=-2L(R-\epsilon),
\end{equation}
and, if $J(R+\epsilon)>0$ for $\epsilon>0$, then
\begin{equation}\label{4.23}
\lim_{N\to\infty}\frac 1{N}\log P_{N,m}^{(\alpha,\beta)}[\frac
1m\max h_j\ge R+\epsilon]=-2J(R+\epsilon).
\end{equation}
We always have $L(R-\epsilon)>0$ if $\epsilon>0$, but we must prove
that $J(R+\epsilon)>0$ if $\epsilon>0$. The function $J$ is defined by
\begin{equation}
J(x)=\inf_{\tau\ge x} g(x),
\notag
\end{equation}
where
\begin{equation}
g(x)=\int_0^\infty\log |x-y|^{-1}u_{\text{eq}}(y)dy+\frac 12V(x)+\frac
12\int_0^\infty V(y) u_{\text{eq}}(y)dy -F_V.
\end{equation}
By the general theory for the constrained variational problem
(\ref{4.21}), see \cite{DS1}, $g(x)\ge 0$ for $x\ge R$. Now, for $x>R$,
\begin{equation}
g''(x)=\int_0^\infty\frac{u_{\text{eq}}(y)}{(x-y)^2}dy+\frac 12V''(x)
\notag
\end{equation}
and $V''(x)=-tU''(tx)$ with
\begin{equation}
U''(s)=-\frac{\alpha_0}{(1+\alpha_0-s)(1-s)}-\frac{\beta_0}{(\beta_0+s)s},
\notag
\end{equation}
so $V''(x)>0$ and $g$ is strictly convex. Consequently $g(x)>0$ if
$x>R$, which is what we wanted to prove.

Note the asymmetry in (\ref{4.22}) and (\ref{4.23}). Just as for the
arctic circle large inward fluctuations of the temperate zone have
much smaller probability than large outward fluctuations.

We will now compute (a part of) the arctic ellipse. For simplicity we
restrict to the case $a=b$. The general case can be handled similarly
but the computations are somewhat more involved. One approach is to
compute the equilibrium measure, and hence its suport which gives the
arctic ellipse, by solving the variational problem (\ref{4.21}) as was
done for the Krawtchouk polynomials in \cite{DS2}. Below we will
instead use the approach of \cite{KvA} which is based on the recursion
formula for the Hahn polynomials.
\begin{lemma}\label{lem4.2}
Let $R(t,\alpha_0)$ be the right endpoint of the support of
$u_{\text{eq}}^{(t,\alpha_0,\alpha_0)}$ as defined above. Then,
\begin{equation}\label{4.24}
R(t,\alpha_0)=\frac 1t\sup_{0<s<t}(\frac 12+\frac
1{2(s+\alpha_0)}\sqrt{ s(1-s)(s+2\alpha_0)(s+2\alpha_0+1)}).
\end{equation}
\end{lemma}
\begin{Proof}
Let $q_n=q_{n,N}^{(\alpha,\beta)}(x)$ denote the normalized Hahn
polynomials which are orthonormal on $\{0,\dots,N\}$ with respect to
the weight (\ref{4.20'}),
\begin{align}\label{4.21'}
q_{n,N}^{(\alpha,\beta)}(x)&=\frac{(-1)^n}{d_{n,N}}{}_3F_2(-n,-x,
n+\alpha+\beta+1;\beta+1,-N;1)\notag\\
&=\frac{(-1)^n}{d_{n,N}}\sum_{k=0}^n\frac{(-n)_k(-x)_k(n+\alpha+\beta+1)_k}{(\beta+1)_k(-N)_kk!},
\end{align}
where $d_{n,N}>0$ and
\begin{equation}
d_{n,N}^2=\frac{(\alpha+\beta+1)(\alpha+1)_n(N+\alpha+\beta+2)_n}
{\binom{N}{n}(2n+\alpha+\beta+1)(\beta+1)_n(\alpha+\beta+1)_n}.
\notag
\end{equation}
The leading coefficients are
\begin{equation}\label{4.25}
\kappa_{n,N}=\frac{(-1)^n}{d_{n,N}}
\frac{(n+\alpha+\beta+1)_n}{(\alpha+1)_n(-N)_n}.
\end{equation}
The polynomials $q_n$ satisfy the recurrence relation
\begin{equation}
xq_n=a_{n,N}q_{n-1}+b_{n,N}q_n+a_{n+1,N}q_{n+1},
\notag
\end{equation}
where
\begin{equation}
a_{n,N}=\frac{n(n+\alpha)(n+\alpha+\beta+N+1)} 
{(2n+\alpha+\beta)(2n+\alpha+\beta+1)}\sqrt{\frac{(N-n+1)(2n+\alpha+\beta+1)
    (\beta+n)(\alpha+\beta+n)}{(\alpha+n)(n+N+\alpha+\beta+1)n
(2n+\alpha+\beta+1)}}
\notag
\end{equation}
and
\begin{equation}
b_{n,N}=\frac{(n+\alpha+\beta+1)(n+\beta+1)(N-n)}{N(2n+\alpha+\beta+1)
  (2n+\alpha+\beta+2)}
  +\frac{n(n+\alpha)(n+\alpha+\beta+N+1)}
{N(2n+\alpha+\beta)(2n+\alpha+\beta+1)} .
\notag
\end{equation}
Consider the rescaled equilibrium measure,
$u_{\text{eq}}^\ast(s)=\frac 1t u_{\text{eq}}(s/t)$, $0\le s\le 1$. The
paper \cite{KvA} tells us how to compute the support of
$u_{\text{eq}}^\ast$ (and also the measure itself) using the
asymptotics of the recursion coefficients. We restrict to the case
$\alpha_0=\beta_0$. Then, $n/N\to t\in (0,1)$,
\begin{align}
\lim_{N\to\infty} a_{n,N}&=a(t)=\frac 1{4(t+\alpha_0)}
\sqrt{t(1-t)(t+2\alpha_0)(t+2\alpha_0+1)}\notag\\
\lim_{N\to\infty} b_{n,N}&=b(t)=\frac 12.
\notag
\end{align}
According to \cite{KvA}, p. 171, the right endpoint of the support of 
$u_{\text{eq}}^\ast$, i.e. $tR(\alpha_0,t)$ is 
\begin{equation}
tR(\alpha_0,t)=\sup_{0<s<t}(b(s)+2a(s)).
\notag
\end{equation}
\end{Proof}

We can now formulate the result we obtain for the arctic ellipse.
\begin{theorem}\label{thm4.3}
Consider the abc-hexagon, assume that $a=b$ and rescale the size of
the hexagon by a factor $1/c$. We let the size of the hexagon grow in
such a way that $a/c\to\lambda>0$ as $c\to\infty$. Pick the
ON-coordinate system in which the limiting hexagon has corners at $\pm
(-\frac{\sqrt{3}}2\lambda,\frac 12)$, $\pm(0,\frac 12(1+\lambda))$ and $\pm(
\frac{\sqrt{3}}2\lambda,\frac 12)$. Let $X$ be the intersection point
of the line $x=\tau$, $-\frac{\sqrt{3}}2<\tau<-\frac{\sqrt{3}}2
\frac{\lambda^2}{\lambda+1}$, with the inner part of the boundary of
the rescaled polar zone at $P_6$. Then,
\begin{equation}\label{4.26}
X\to\sqrt{2\lambda+1}\sqrt{1/4-\tau^2/3\lambda^2}
\end{equation}
in probability as $c\to\infty$. For any $\epsilon>0$ there are
constants $I(\epsilon)>0$ and $J(\epsilon)>0$ such that
\begin{equation}\label{4.27}
\frac 1c\log
P[X\ge\sqrt{2\lambda+1}\sqrt{1/4-\tau^2/3\lambda^2}+\epsilon]\to-I(\epsilon)
\end{equation}
and
\begin{equation}\label{4.28}
\frac 1{c^2}\log
P[X\le\sqrt{2\lambda+1}\sqrt{1/4-\tau^2/3\lambda^2}-\epsilon]\to-J(\epsilon)
\end{equation}
as $c\to\infty$.
\end{theorem}
\begin{Proof}
We are in the case when $0\le m\le a=b$, so $n=L_m=m$,
$N=\gamma_m=m+c-1$, $\alpha=a-m=\beta=b-m$. Assume that $m/c\to\mu>0$
as $c\to\infty$. Then $t=\mu(1+\mu)^{-1}$ and
$\alpha_0=(\lambda-\mu)(1+\mu)^{-1}=(1-t)\lambda-t$. Consider the
points $B_m=(m,\alpha_m+2Z_m)$ which describe the inner boundary of
the polar zone; $\alpha_m=-m$. We see that
\begin{equation}\label{4.29}
Z_m/L_m\to R(t,\alpha_0)
\end{equation}
in probability as $c\to\infty$ Here we use the large deviation
formulas (\ref{4.22}) and (\ref{4.23}) together with \ref{thm4.1}. To
get an ON-system (with the sides of the rhombi $=1$ we have to rescale
the coordinates to $\tilde{B}_m=(\frac{\sqrt{3}}2m,-\frac
m2+Z_m)$. Thus,
\begin{equation}
\frac 1c\tilde{B}_m=(\frac{\sqrt{3}}2\frac mc,-\frac
m{2c}+\frac{L_m}c\frac{Z_m}{L_m}) \to
(\frac{\sqrt{3}}2\mu,-\frac{\mu}2+\mu R)
\notag
\end{equation}
We also have to translate the coordinate system so that we get the
origin at the center of the hexagon. We then get the coordinates
\begin{equation}\label{4.30}
(\frac{\sqrt{3}}2(\mu-\lambda),-\frac{\mu+1}2+\mu R)
\end{equation}
for the limit of the point $X$ on the arctic ellipse. Now,
\begin{equation}\label{4.31}
R(\frac{\mu}{\mu+1},\frac{\lambda-\mu}{\mu+1})=\frac{\mu+1}{2\mu}
+\frac 1{2\lambda\mu}\sqrt{(2\lambda+1)\mu(2\lambda-\mu)}
\end{equation}
for $0\le\mu\le\lambda(\lambda+1)^{-1}$. To see this we use lemma
\ref{lem4.2} and set
\begin{equation}
g(s)=\frac 12+\frac
1{2(s+\alpha_0)}\sqrt{s(1-s)(s+2\alpha_0)(s+2\alpha_0+1)}
\notag
\end{equation}
with $\alpha_0=(1-t)\lambda-t$. We must have $g(s)\le 1$ because the
support is restricted to $[0,1/t]$. A computation shows that $g(s)$ is
strictly increasing in $[0,s_0]$ and strictly decreasing in $[s_0,1]$
for some $s_0\in(0,1)$ For $0\le t\le s_0$,
\begin{equation}\label{4.32}
\sup_{0<s<t} g(s)=g(t)=\frac 12+\frac{2\lambda+1}{2\lambda}
\sqrt{t(\frac{2\lambda}{2\lambda+1}-t)}, 
\end{equation}
where we have inserted $\alpha_0=(1-t)\lambda-t$. Suppose
$s_0<\lambda(2\lambda+1)^{-1}$. Then
$g(s_0)>g(\lambda(2\lambda+1)^{-1})=1$, which is impossible. Thus
(\ref{4.32}) holds for $0\le t\le \lambda(2\lambda+1)^{-1}$ and
(\ref{4.31}) follows from (\ref{4.24}) since $t=\mu(\mu+1)^{-1}$. Set $\tau=
\frac{\sqrt{3}}2(\mu-\lambda)\in(-\frac{\sqrt{3}}2,-\frac{\sqrt{3}}2
\frac{\lambda^2}{\lambda+1})$. Then the lmiting point (\ref{4.30})
becomes $(\tau,\sqrt{2\lambda+1}\sqrt{1/4-\tau^2/3\lambda^2})$. The
large deviation formulas (\ref{4.27}) and (\ref{4.28}) follow from
(\ref{4.22}) and (\ref{4.23}).
\end{Proof}

The 1-dimensional marginal probability in the Hahn ensemble
(\ref{4.3}) is, \cite{Me},
\begin{equation}\label{4.33}
u_{N,n}^{(\alpha,\beta)}(t)=\frac
1n\sum_{k=0}^{n-1}q_{k,N}^{(\alpha,\beta)} (t)^2W_N^{(\alpha,\beta)}(t),
\end{equation}
for $t\in\{0,\dots,N\}$, where $q_{k,N}^{(\alpha,\beta)}$ are the
normalized Hahn polynomials (\ref{4.21'}). Hence, the probability of
finding a hole at position $t$, i.e. the 1-point correlation function,
is $nu_{N,n}^{(\alpha,\beta)}(t)$, and consequently the number of
rhombus tilings of the abc-hexagon with a vertical rhombus at position
$t$ on the line $x=m$ is
\begin{equation}\label{4.34}
N(a,b,c)\sum_{k=0}^{L_m-1}q_{k,N}^{(\alpha,\beta)}
(t)^2W_N^{(a_m,b_m)}(t) ,
\end{equation}
by theorem \ref{thm4.1}. (This can be rewritten using the
Christoffel-Darboux formula.) If we use the explicit formula
(\ref{4.21'}) for the Hahn polynomials, we obtain an explicit formula
for the quantity (\ref{4.34}). This quantity has been investigated in
\cite{Fisch}, \cite{FK1}, \cite{FK2}. 
We will not attempt to rewrite (\ref{4.34}) in oeder to compare
it with existing formulas. We can also consider the number of tilings
with vertical rhombi at specified positions $t_1,\dots,t_r$ on the
line $x=m$. This is given by $N(a,b,c)$ times a
determinantal correlation function like (\ref{2.8}) but where we now
have instead the Hahn kernel given by the formula (\ref{2.9}) with the
Hahn polynomials and the Hahn weight instead. It follows from the
general theory in \cite{Jo1} that the 1-dimensional marginal
probability converges weakly to the equilibrium measure. This should
also hold pointwise, i.e.
\begin{equation}\label{4.35}
\lim_{N\to\infty}u_{N,n}^{(\alpha,\beta)}([n\tau])
=u_{\text{eq}}^{(t,\alpha_0,\beta_0)}(\tau),
\end{equation}
if $\frac 1N (\alpha,\beta)\to(\alpha_0,\beta_0)$ and $n/N\to t$ as
$N\to\infty$, but we do not have good enough control over the
asymptotics of the Hahn polynomials to prove it. The corresponding
result for the Krawtchouk ensemble follows by the same methods as was
used to prove lemma \ref{lem2.8}.

As mentioned in the beginning of this section a rhombus tiling of an
abc-hexagon can also be interpreted as a boxed planar partition where
the sides of the box are $a$, $b$ and $c$. We want to relate the
height of the planar partition surface above a certain point to our
particle configurations. Let $(x,y,z)$ be the coordinates in the
planar partition coordinate system (the $x$-axis through $P_1$, the
$y$-axis through $P_3$ and the $z$-axis through $P_5$ after
projection). Consider the line $x=m$. We start from the point
$(r_m,s_m,0)$, where
\begin{equation}
r_m=
\begin{cases}
a&\text{if $0\le m\le b$}\\
a+b-m&\text{if $b\le m\le a+b$}
\end{cases}
\notag
\end{equation}
and
\begin{equation}
s_m=
\begin{cases}
m&\text{if $0\le m\le b$}\\
b&\text{if $b\le m\le a+b$}
\end{cases}.
\notag
\end{equation}
As we go along the line a particle means that the $z$-coordinate is
increased by 1, whereas the $x$- and $y$-coordinates are fixed. A hole
means that the $z$-coordinate is fixed, but the $x$- and
$y$-coordinates are reduced by 1. Let $X_m(k)$ denote the position of
hole number $k$. This hole corresponds to the position $(r_m-k,s_m-k)$
in the $xy$-plane and the surface height above this point is equal to
the number of of particles in $\{0,\dots,X_m(k)\}$, which equals
$X_m(k)-k+1$. Thus, the planar partition height function is given by
\begin{equation}\label{4.36}
H(r_m-k,s_m-k)=X_m(k)-k+1.
\end{equation}
Let $Y_m(n)$ denote the number of holes in $[0,n]$ in the particle
system on $x=m$. Then
\begin{equation}\label{4.37}
P[X_m(k)\le n]=1-P[Y_m(n)<k].
\end{equation}
Assume that $m/c\to\mu >0$, $L_m/\gamma_m\to t>0$, $n/L_m\to\tau$ and
$\gamma^{-1}(|a-m|,|b-m|)\to(\alpha_0,\beta_0)$, as $c\to\infty$. Then
theorem \ref{thm4.1} and the general theory of discrete Coulomb gases
in \cite{Jo1} shows that
\begin{equation}\label{4.38}
\frac
1{L_m}E[Y_m(n)]\to\int_0^{\tau}u_{\text{eq}}^{(t,\alpha_0,\beta_0)}(s)ds 
\end{equation}
as $c\to\infty$. Also, it is possible to prove large deviation
results, analogous to those in \cite{BG} for the semi-circle law, in
this case too. Using this and (\ref{4.36}), (\ref{4.37}) and
(\ref{4.38}) it is possible to compue the asymptotic shape of the
planar partition surface (law of large numbers) as well as large
deviation results. See \cite{CLP} for this type of results proved by other
methods. The asymptotic shape can be computed if we know the
equilibrium measure.

From (\ref{4.36}) and (\ref{4.37}) it is clear that the surface
fluctuations are directly related to the fluctuations of the number of
paticles in an interval in the Hahn ensemble. Hence, in analogy with
the Krawtchouk case, the surface fluctuations should be Gaussian with
variance proportional to $\log c$ as $c\to\infty$. This could be
proved provided we had the same control of the Hahn kernel as we have
of the Krawtchouk kernel in lemma \ref{lem2.8}, but this remains to be
done.

We will end this section by showing one way of finding GUE in a random
rhombus tiling of an abc-hexagon in the limit as the size of the
hexagon goes to infinity. Consider the $m$ holes in the $m$:th column,
$1\le m\le b= a$. we want to compute the probability distribution of
these $m$ holes, with $m$ fixed, as the hexagon grows. Let
$\xi_1,\dots,\xi_m$ be the positions of the holes. By theorem
\ref{thm4.1} the probability of this configuration is
\begin{equation}\label{4.39}
\frac 1{Z_{\gamma_m,m}^{(a-m,a-m)}}\Delta_m(\xi)^2\prod_{j=1}^m
\frac{(\gamma_m+a-m-\xi_j)!(a-m+\xi_j)!}{\xi_j!(\gamma_m-\xi_j)!},
\end{equation}
where $\gamma_m=m+c-1$. Denote the corresponding expectation by
$E_{\gamma_m,m}^{(a-m,a-m)}[\cdot]$.
\begin{proposition}\label{prop4.4}
Let $c\to\infty$, $a/c\to\lambda>0$ and keep $m\ge 1$ fixed. Let
$f:\mathbb{R}^m\to\mathbb{C}$ be a continuous, bounded, symmetric
function. Then,
\begin{align}\label{4.40}
&\lim_{c\to\infty}E_{\gamma_m,m}^{(a-m,a-m)}[f(\frac{\xi_1-\gamma_m/2}{\sqrt{c}},
\dots, \frac{\xi_m-\gamma_m/2}{\sqrt{c}})] \notag\\
&=\frac 1{Z_m(\lambda)}\int_{\mathbb{R}^m}\Delta_m(x)^2\prod_{j=1}^m
e^{-\frac{2\lambda}{\lambda+1}x_j^2} f(x_1,\dots,x_m)d^mx,
\end{align}
where $Z_m(\lambda)$ is a normalization constant such that the right
hand side is 1 when $f=1$.
\end{proposition}

Thus, in the limit, the positions of the holes (vertical rhombi) on
the $m$:th vertical column are described by $m\times m$ GUE.

\begin{Proof}
The proof is straightforward. By (\ref{4.39}) the expectation in the
left hand side of (\ref{4.40}) equals
\begin{align}
\frac 1{Z_{\gamma_m,m}^{(a-m,a-m)}}
\sum_{\xi\in\{0,\dots,\gamma_m\}^m} &\Delta_m(\xi)^2
f(\frac{\xi_1-\gamma/2}{\sqrt{c}},
\dots, \frac{\xi_m-\gamma/2}{\sqrt{c}})\notag\\
&\times\prod_{j=1}^m
\frac{(\gamma+a-m-\xi_j)!(a-m+\xi_j)!}{\xi_j!(\gamma_m-\xi_j)!}.
\notag
\end{align}
In this expression we write $\xi_j=n_j+\gamma/2$ and use Stirling's
formula to approximate the factorials. This leads to an expression
which is a Riemann sum. The normalization constant is the same Riemann
sum but with $f=1$. After cancelling common factors we see that the
remaining quotient of Riemann sums converges to the right hand side of
(\ref{4.40}).
\end{Proof}

\section{A dimer model on a cylindrical brick lattice}

In this final section we will consider a dimer model on a hexagonal
graph on a finite cylinder. The graph is sometimes referred to as the
brick lattice, \cite{YNS}. A dimer covering of this lattice can also
be thought of as a certain cylindrical rhombus tiling, see
\cite{Ke1}. The dimer covering has an equivalent description in terms
of non-intersecting random walk paths, \cite{Fi}. In contrast to the
previous non-intersecting path models we do not have a fixed number of
paths. We will again be interested in the point process we obtain by
looking at where the random walks are at a fixed time, and the
analysis is based on the methods of \cite{TW2}. 

The graph $G_{M,N}$ which we will consider is defined as follows. The
vertices are $v_{j,k}=(-1/2+j,k)$, $0\le j\le 2M-1$, $0\le k\le 2N$
and we obtain a graph on a cylinder by identifying $v_{j,k}$ and
$v_{j+2m,k}$ for all $j,k$. We have vertical edges between $v_{j,k}$
and $v_{j,k+1}$, and horizontal edges between $v_{2j,2k}$ and
$v_{2j+1,2k}$ and between $v_{2j+1,2k+1}$ and $v_{2j+2,2k+1}$. A dimer
covering of $G_{M,N}$ can equivalently be described by
non-intersecting random walk paths, \cite{Fi},\cite{NYB}. In the same
coordinate system as was used to define $G_{M,N}$ we have simple
random walk paths $S(t)$ with steps $\pm 1$, with initial position
$S(0)\in \{2k\,;\, 0\le k\le N\}$ and which are required to satisfy
$0\le S(t)\le 2N$. Since we have a graph on a cylinder we and the
paths to live on the cylinder also so we require that
$S(t+2M)=S(t)$. We have $L$ non-intersecting paths, where $0\le L\le
N$. (Actually the condition $0\le S(t)\le 2N$ implies that $L\le
N-1$.) Let $\mathcal{P}_{M,N}$ be th set of all such families of
non-intersecting paths. We now describe the 1-1-correspondence between
$\mathcal{P}_{M,N}$ and the set of all dimer configurations on
$G_{M,N}$. A dimer covers a vertical edge if and only if we have a
random walk step, in one of the paths, which intersects the vertical
edge. A dimer covers a horizontal edge if and only if {\it no} random
walk path hits this edge. If the horizontal edge
$v_{2j,2k}v_{2j+1,2k}$ is {\it not} covered by a dimer, then one of
$e_1=v_{2j+1,2k}v_{2j+1,2k+1}$ or $e_2=v_{2j+1,2k}v_{2j+1,2k-1}$ must
be covered. A random walk path must pass through $(2j,2k)$ and in the
next step intersect either $e_1$ or $e_2$. It will go to either
$(2j+1,2k+1)$ or $(2j+1,2k-1)$ which means that we obtain a new
horizontal edge that is not covered by a dimer and we can repeat the
argument. A random walk path connot inresect both $e_1$ and $e_2$
because then they would have to meet at $(2j,2k)$ which would
contradict the non-intersection condition. A similar argument applies
at all locations in the graph. Hence, a dimer configuration gives rise
to non-intersecting paths, and conversely if we have non-intersecting
paths we can find the dimer covering corresponding to them.

Consider a dimer covering of $G_{M,N}$. The total number of vertical 
dimers is $2ML$, where $L$ is the number of non-intersecting paths in
the path description. The total number of dimers is $M(2N+1)$ and
hence the numbeer of horizontal dimers is $M(2N+1-2L)$. We will now
define a probability measure on the set of dimer configurations by
letting vertical dimers have weight $w$ and horizontal dimers weight
$z$. Let $g_{M,N}(m,n)$ denote the number of configurations with $m$
horizontal and $m$ vertical dimers; each of these have the same
probability. The {\it partition function} is given by
\begin{equation}\label{5.1}
Z=Z_{M,N}(z,w)=\sum_{m,n\ge 0} g_{M,N}(m,n)z^mw^n.
\end{equation}
If $G_L$ is the number of non-intersecting path configurations with
exactly $L$ paths, then we must also have
\begin{equation}\label{5.2}
Z=\sum_{L=0}^NG_Lz^{M(2N+1-2L)}w^{2ML}=z^{M(2N+1)}\sum_{L=0}^NG_L
\left(\frac wz\right)^{2ML}.
\end{equation}
Recall that the possible initial (=final) positions for the
non-intersecting paths are $\{2k\,;\, 0\le k\le N\}$. Let $G_L(x)$
denote the number of configurations with $L$ non-intersecting paths
whose initial conditions are $2x_1<\dots<2x_L$, where
$x_1,\dots,x_L\in [N]\doteq \{0,\dots,N\}$ are given. Then,
\begin{equation}\label{5.3}
G_L=\sum_{0\le x_1<\dots<x_L\le N} G_L(x).
\end{equation}
Hence, the probability that the initial positions $x$
{\it given} that we have exactly $L$ non-intersecting paths is
$G_L(x)/G_L$, and we define a probability on $[N]^L$ by
\begin{equation}\label{5.4}
u_L(x)=\frac 1{L!G_L}G_L(x),
\end{equation}
where the right hand side is extended to $[N]^L$ by requiring it to be
a symmetric function. The {\it $l$-particle correlation function},
{\it given} that the total number of particles is $L$, is defined by
\begin{equation}\label{5.5}
R_{\ell,L}(x_1,\dots,x_\ell)=\frac{L!}{(L-\ell)!}
\sum_{x_{\ell+1},\dots,x_L\in[N]} u_L(x_1,\dots,x_L).
\end{equation}
The probability of having exactly $L$ particles is, by (\ref{5.2})
\begin{equation}
\frac {G_L}{Z}(w/z)^{2ML}z^{M(2N+1)}\notag
\end{equation}
and hence the {\it $\ell$-particle correlation function}, with no
restriction on the total number of particles, is
\begin{equation}\label{5.6}
R_\ell(x_1,\dots,x_\ell)=\frac{z^{M(2N+1)}}{Z}\sum_{L=\ell}^N
R_{\ell,L}(x_1,\dots,x_\ell)(w/z)^{2ML}G_L.
\end{equation}
Set
\begin{equation}\label{5.7}
\phi(s,t)=\sqrt{\frac{2c_{2s}c_{2t}}{N}}\cos\frac{\pi st}N,
\end{equation}
$0\le s,t\le N$, where $c_m=1/2$ if $m=0$ or $m=2N$ and $c_m=1$ if
$1\le m\le 2N-1$.

\begin{proposition}\label{prop5.1}
Set
\begin{equation}\label{5.8}
K(x,y)=\sum_{j=0}^N\phi(x,j)\phi(y,j)\frac{(2w/z)^{2M}w_j}{1+(2w/z)^{2M}w_j},
\end{equation}
where $w_j=(\cos\frac{\pi j}{2N})^{2M}$. Then,
\begin{equation}\label{5.9}
R_\ell(x_1,\dots,x_\ell)=\det(K(x_i,x_j))_{i,j=1}^\ell.
\end{equation}
Also,
\begin{equation}\label{5.9'}
Z_{M,N}(z,w)=z^{M(2N+1)}\prod_{k=0}^N
(1+(\frac{2w}z\cos\frac{\pi k}{2N})^{2M}).
\end{equation}
\end{proposition}
\begin{Proof}
We consider a simple random walk on $\{-1,0,\dots,2N+1\}$, where $-1$
and $2N+1$ are absorbing barriers. The transition matrix is given by
$P^\ast_{n,n}=0$ if $0\le n\le 2N$, $P^\ast_{n,n}=1$ if $n=-1$ or
$n=2N+1$, $P^\ast_{n,n-1}=1/2$ if $1\le n\le 2N+1$
$P^\ast_{n,n+1}=1/2$ if $0\le n\le 2N-1$ and $P^\ast_{m,n}=0$
otherwise. Since we are only interested in random walk paths that are
not absorbed we can concentrate on the
submatrix $P=(P^\ast_{m,n})_{0\le m,n\le 2N}$. We are interested in the
probability of going from $2x$ to $2y$ in $2M$ steps, which is given
by $(P^{2M})_{2x,2y}$. A computation, see \cite{Karl1}, shows that
\begin{equation}
\sum_{j=0}^N(P^{2M})_{2x,2j}\phi(j,k)=w_k\phi(x,k)
\notag
\end{equation}
and because of the orthogonality
\begin{equation}
\sum_{j=0}^N\phi(s,j)\phi(j,t)=\delta_{st},
\notag
\end{equation}
we obtain
\begin{equation}\label{5.10}
(P^{2M})_{2x,2y}=\sum_{j=0}^Nw_j\phi(x,j)\phi(j,y),
\end{equation}
with $w_j$ as above.

Now, it follows from the Karlin/McGregor theorem that
\begin{equation}\label{5.11}
G_L(x)=2^{2ML}\det((P^{2M})_{2x_i,2x_j})_{i,j=1}^L,
\end{equation}
where the factor $2^{2ML}$ comes from the fact that $G_L(x)$ counts
the {\it number} of configurations. Note that the right hand side of
(\ref{5.11}) is a symmetric function of $x_1,\dots,x_L$. Let $f$ be a
given function on $\mathbb{N}$ and set
\begin{equation}\label{5.12}
\Phi_L[f]=\sum_{x\in
  [N]^L}\det((P^{2M})_{2x_i,2x_j})_{i,j=1}^L\prod_{j=1}^L f(x_j),
\end{equation}
and write $[\lambda^\ell]F(\lambda)$ for the coefficient of
$\lambda^\ell$ in the power series $F(\lambda)$. It follows from (\ref{5.4}),
(\ref{5.5}), (\ref{5.11}) and (\ref{5.12}) that
\begin{equation}\label{5.13}
\sum_{x\in [N]^L}R_{\ell,L}(x_1,\dots,x_{\ell})\prod_{j=1}^L f(x_j)
=\frac{\ell!2^{2ML}}{L!G_L}[\lambda^\ell]\Phi_L[1+\lambda f].
\end{equation}
Also, we set
\begin{equation}\label{5.14}
\Phi[f]=\sum_{L=0}^N\frac{(2w/z)^{2ML}}{L!}\Phi_L[f]
\end{equation}
and note that by (\ref{5.2}), (\ref{5.3}), (\ref{5.11}) and
(\ref{5.12})
\begin{equation}\label{5.15}
\Phi[1]=\frac {Z}{z^{M(2N+1)}}.
\end{equation}
Now, by (\ref{5.6}), (\ref{5.13}), (\ref{5.14}) and (\ref{5.15}),
\begin{equation}\label{5.16}
\sum_{x\in[N]^\ell}R_\ell(x_1,\dots,x_\ell)
\prod_{j=1}^\ell f(x_j)=\ell![\lambda^\ell]
\frac{\Phi[1+\lambda f]}{\Phi[1]}.
\end{equation}
By (\ref{5.10}) and a classical identity, see \cite{TW2},
\begin{align}
\det((P^{2M})_{2x_i,2x_j})_{i,j=1}^L&=
\det(\sum_{k=0}^Nw_k\phi(x_i,k)\phi(k,x_j))_{i,j=1}^L\notag\\
&=\frac 1{L!}\sum_{k\in[N]^L}[\det(\phi(x_i,k_j)_{i,j=1}^L]^2
\prod_{j=1}^L w_{k_j}.\notag
\end{align}
Thus, if we use the same identity again in the other direction we obtain
\begin{equation}\label{5.17}
\Phi_L[f]=\sum_{k\in[N]^L}\det(\sum_{x=0}^N w_{k_i}\phi(x,k_i)\phi(x,k_j) 
f(x))_{i,j=1}^L.
\end{equation}
Set
\begin{equation}
\mathcal{K}(s,t)=\sum_{x=0}^N w_{s}\phi(x,s)\phi(x,t) 
f(x).\notag
\end{equation}
Then, by (\ref{5.14}), (\ref{5.17}) and a Fredholm expansion
\begin{equation}\label{5.18}
\Phi[f]=\det(\delta_{i,j}+(2w/z)^{2M}\mathcal{K}(i,j))_{i,j=0}^N.
\end{equation}
Set $A=((1+(2w/z)^{2M}w_i)\delta_{i,j})_{i,j=0}^N$ and 
$B=((2w/z)^{2M}\mathcal{K}(i,j))_{i,j=0}^N$. Then, by (\ref{5.18}),
\begin{equation}\label{5.19}
\Phi[1+\lambda f]=\det(A+\lambda B).
\end{equation}
Combining (\ref{5.15}) and (\ref{5.19}) we obtain (\ref{5.9'}). 
Furthermore, (\ref{5.16}) can be written
\begin{equation}\label{5.20}
\sum_{x\in[N]^\ell}R_\ell(x_1,\dots,x_\ell)\prod_{j=1}^\ell f(x_j)
=\ell![\lambda^\ell]\det(I+\lambda A^{-1}B)
=\ell![\lambda^\ell]\det (I+\lambda Kf)
\end{equation}
with $K$ given by (\ref{5.8}). The last step is proved as in \cite{TW2}, 
see also \cite{Jo3}. Expanding the last expression in (\ref{5.20}) in a 
Fredholm expansion gives (\ref{5.9}) since $f$ was arbitrary, and the 
proposition is proved 
\end{Proof}

We will now discuss the asymptotics of the partition function and the 
correlation functions. It is well known how to obtain a formula like 
(\ref{5.9'}) for the partition function using Kasteleyn's method, \cite{K}.
Also the two-point correlation function has been computed in \cite{YNS} 
using the methods of \cite{FiSt}. See also \cite{Fo2'}. The computation
above is different and emphasizes the similarity between certain aspects 
of the dimer model and random matrix theory. The limiting correlation
functions we obtain are given by the discrete or ordinary sine
kernel. 

The number of vertices is $2M(2N+1)$ and the {\it free energy} per
vertex is defined by
\begin{equation}\label{5.21}
f_{M,N}(z,w)=\frac 1{2M(2N+1)}\log Z_{M,N}(z,w).
\end{equation}
\begin{proposition}\label{prop5.2}
The limiting free energy is
\begin{align}\label{5.22}
f(z,w)&=\lim_{M\to\infty}\lim_{N\to\infty}f_{M,N}(z,w)\notag\\
&=
\begin{cases} 
\frac 12\log z & \text{if $w/z<1/2$}\\
\frac 12\log z+\frac 1{2\pi}\int_0^{\pi\theta_0/2}
\log(\frac{2w}z\cos s)ds  & \text{if $w/z>1/2$},
\end{cases}
\end{align}
\end{proposition}
\begin{Proof}
By (\ref{5.9'}) and (\ref{5.21}) we obtain
\begin{align}
\lim_{N\to\infty}f_{M,N}(z,w)&=\frac 12\log z+\frac
1{4M}\lim_{N\to\infty} \frac 1{N+1}\sum_{k=0}^N\log
(1+(\frac{2w}z\cos\frac{\pi k}{2N})^{2M}) \notag\\
&=\frac 12\log z+\frac 1{2\pi M}+\int_0^{\pi /2} \log(1+(\frac{2w}z\cos
s)^{2M}) ds\notag
\end{align}
If $w/z<1/2$, then the integrand goes to 0 uniformly and we obtain the
first part of (\ref{5.22}). If $w/z>1/2$, then the integrand goes to 0
unless $0\le s\le \arccos (z/2w)=\pi\theta_0/2$ and we obtain our
result by using the inequalities
\begin{equation}
(\frac{2w}z\cos s)^{2M}\le 1+(\frac{2w}z\cos s)^{2M} \le
2(\frac{2w}z\cos s)^{2M} 
\end{equation}
for $0\le s\le \theta_0$.
\end{Proof}

If we only had horizontal dimers ($L=0$), then $Z=z^{M(2N+1)}$ and
hence the free energy per vertex is $\frac 12\log z$. We can thus
interpret the phase transition in (\ref{5.22}) as saying that for
$w/z<1/2$ the system is completely frozen, whereas for $w/z>1/2$ we
have both horizontal and vertical dimers. This type of phase
transition is called a {\it K-type transition} in \cite{NYB}. We will study
the limiting correlation functions in the non-frozen phase $w/z>1/2$
and in a scaling limit where we approach the critical point from
above. In this scaling limit we will obtain the sine kernel
determinantal point process of random matrix theory. By proposition
\ref{prop5.1}, formula (\ref{5.9}) it suffices to investigate the
asymptotics of the kernel (\ref{5.8}). 

\begin{proposition}\label{prop5.3}
Assume that $w/z>1/2$. Then
\begin{equation}\label{5.23}
\lim_{M\to\infty}\lim_{N\to\infty} K(\frac N2+t,\frac N2+s)=\frac{
\sin\pi (t-s)\theta_0}{\pi (t-s)},
\end{equation}
for any fixed $t,s\in\mathbb{Z}$. Here $\theta_0=\frac 2{\pi}\arccos
\frac z{2w}$ as before. Let $\epsilon_N$, $N\ge 1$, be a given
sequence such that $\epsilon\to 0$ and $N\sqrt{\epsilon_N}\to \infty$
as $N\to\infty$. Assume furthermore that $M=M(N)\to\infty$ and
$M\epsilon_N\to\infty$ as $N\to\infty$. If $2w/z=1+\epsilon_N$, then
\begin{equation}\label{5.24}
\lim_{N\to\infty} \frac{\pi}{2\sqrt{2\epsilon_N}}K(\frac N2+[\frac
{\pi\xi}{2\sqrt{2\epsilon_N}}],\frac N2+[\frac
{\pi\eta}{2\sqrt{2\epsilon_N}}])=\frac{ 
\sin\pi (\xi-\eta)}{\pi (\xi-\eta)},
\end{equation}
for any fixed $\xi,\eta\in\mathbb{R}$.
\end{proposition}

Note that $\theta_0$ gives the local density of the intersection of
the random walk paths with a vertical axis and that $\theta\to 0$ as
$w/z\to1/2+$.

\begin{Proof}
Assume that $N=2n$ for simplicity. Then, by the definition of
$K(x,y)$,
\begin{align}
&\lim_{N\to\infty} K(n+t,n+s)\notag\\
&\lim_{n\to\infty}\frac 1n\sum_{j=0}^{2n}\cos (\frac{\pi j}2+
\frac{\pi tj}{2n})\cos (\frac{\pi j}2+
\frac{\pi sj}{2n})\frac{(\frac{2w}z\cos\frac{\pi j}{4n})^{2M}}
{1+(\frac{2w}z\cos\frac{\pi j}{4n})^{2M}}.
\notag
\end{align}
Split the last sum into two depending on whether $j$ is even or odd
and compute the limits of the Riemann sums obtained. This gives
\begin{equation}
\int_0^1[\sin\pi tu\sin\pi su+\cos\pi tu\cos\pi su]
\frac{(\frac{2w}z\cos\frac{\pi u}{2})^{2M}}
{1+(\frac{2w}z\cos\frac{\pi u}{2})^{2M}}du.\notag
\end{equation}
With $\theta_0$ as defined above we see that the limit of this
expression as $M\to\infty$ is
\begin{equation}
\int_0^{\theta_0}\cos\pi (t-s)udu=\frac{
\sin\pi (t-s)\theta_0}{\pi (t-s)},\notag
\end{equation}
and we have proved (\ref{5.23}). Write
$\gamma_N=1+\epsilon_N=2w/z$. We have that
\begin{align}\label{5.25}
K(\frac N2+t,\frac N2+s)\approx \frac 2N\sum_{j=1}^{N/2}
&\left[\sin(\frac{\pi t(2j-1)}N)\sin(\frac{\pi s(2j-1)}N)
\left(\frac{\gamma_N\cos\frac{\pi s(2j-1)}{2N}}
{1+\gamma_N\cos\frac{\pi s(2j-1)}{2N}}\right)^{2M}
\right.\notag \\
+
&\left.\cos\frac{2\pi tj}N\cos\frac{2\pi sj}N
\left(\frac{\gamma_N\cos\frac{\pi j}N}{1+\gamma_N\cos\frac{\pi j}N}
\right)^{2M}\right],
\end{align}
where the error is negligible for large $N$. Note that
$(1+\epsilon_N)\cos\frac{\pi j}N\ge 1$ if (approximately) $j\le
N\sqrt{2\epsilon_N}/\pi$. Hence the summation in (\ref{5.25}) can be
restricted to $1\le j\le N\sqrt{2\epsilon_N}/\pi$ and in the limit we
are considering the right hand side of (\ref{5.25}) becomes
\begin{equation}
\frac{\pi}{\sqrt{2}}\int_0^{\sqrt{2}/\pi}\cos(\frac{\pi^2}{\sqrt{2}}(\xi
-\eta)u)du =\frac{ 
\sin\pi (\xi-\eta)}{\pi (\xi-\eta)},
\notag
\end{equation}
\end{Proof}

\noindent
{\bf Acknowledgement} I would like to thank S. Karlin for drawing
my attention to his work on coincidence probabilities.

\end{document}